\documentclass{amsart}
\usepackage{amsthm}
\usepackage{amsfonts}
\newtheorem{theorem}{Theorem}
\newtheorem{lemma}[theorem]{Lemma}
\theoremstyle{definition}
\newtheorem{definition}[theorem]{Definition}
\theoremstyle{remark}
\newtheorem{remark}[theorem]{Remark}
\numberwithin{equation}{section}
\theoremstyle{corollary}

\theoremstyle{proposition}
\newtheorem{proposition}[theorem]{Proposition}

\newcommand\B{\mathbb{B}}
\newcommand\N{\mathbb{N}}

\newcommand\C{\mathbb{C}}

\newcommand\de{\partial}

\def\Kl1{\mbox{{\rm K-}}\!\lim_{\hspace{-.4cm}z\rightarrow 1}}

\newcommand\id{{\rm id}}
\newcommand\Aut{{\rm Aut }}
\newcommand\Hol{{\rm Hol }}
\begin{document}

\begin{abstract}
The classical Julia-Wolff-Carath\'{e}odory  Theorem is one of the main tools to study the boundary behavior
of holomorphic self-maps of the unit disc of $\C$.
In this paper we prove a Julia-Wolff-Carath\'{e}odory's type theorem in the case of the polydisc of $\C^n.$
The Busemann functions are used to define a class of  ``generalized horospheres"
for the polydisc and to extend the notion of non-tangential limit.
With these new tools we give a generalization of
the classical Julia's Lemma and of the Lindel\"{o}f Theorem, which the new Julia-Wolff-Carath\'{e}odory Theorem
relies upon.
\end{abstract}
%
%
\title[Busemann Functions and JWC Theorem for Polydiscs]{Busemann Functions and Julia-Wolff-Carath\'{e}odory  Theorem for Polydiscs }
\author[Chiara Frosini]{Chiara Frosini$^\S$}
\thanks{\rm $^\S$ Supported by Progetto MIUR di
Rilevante Interesse Nazionale {\it Propriet\`a geometriche delle
variet\`a reali e complesse}, by G.N.S.A.G.A (gruppo
I.N.D.A.M), and by Progetto FIRB-Dinamica e azioni di gruppi su domini e variet\`a.}
\address{Dipartimento di
Matematica ``U. Dini'', Universit\`a di Firenze, Viale Morgagni
67/A, 50134 Firenze , Italy.} \email{frosini@math.unifi.it}
\date{May 22, 2006}

\maketitle


\section{Introduction}

The Julia-Wolff-Carath\'{e}odory Theorem  and its variants are powerful tools for investigating
the boundary behavior of
holomorphic self-maps of the unit disc $\Delta\subset \C$
(see, {\sl e.g.} \cite{LibroAbate}, \cite{BS},
\cite{Car4}, \cite{Julia2}, \cite{Wolff1}, \cite{Wolff2}).
The importance of this classical theorem (JWC's Theorem for
short)
 in different contexts such as the
study of dynamics, extension of biholomorphisms, composition
operators, semigroups of holomorphic maps, is well known and justifies several generalizations
to higher and infinite dimensions due to various authors. We cite here Rudin
\cite{RUDIN} for the case of the unit ball in $\C^n$, Abate (\cite{ArticoloLindelofAbate}, \cite{ArticoloPseudoConvexAbate}) for
strongly convex and strongly pseudoconvex domains, and Reich and Shoikhet \cite{RS} for the infinite
dimensional case.
The argument used to prove the classical JWC's Theorem inspires its generalizations: let
$f$ be a holomorphic self-map of a domain $D$ and $p\in \de D.$ If
the distance of $f(z)$ from the boundary of the domain, $\hbox{dist}(f(z),\de D),$ is
``comparable" to $\hbox{dist}(z,\de D)$ as $z\to p$ (no matter along
which direction), then $f$ together with its normal derivatives has
a limit at $p$ along some {\sl admissible} directions.
We point out that  to find the class of these admissible directions is one of the main efforts to
state  and prove  all JWC's-type theorems.

In this paper we prove a  JWC's Theorem  in the case of
 the polydisc $\Delta^n$ of $\C^n$ which generalizes the one obtained by Abate
 \cite{ArticoloJWCAbate} (see also Jafari \cite{Jafari}).
In order to achieve this result we first prove
 generalizations of the Julia's Lemma and the Lindel\"{o}f Theorem.
  Our main issue is the use of Busemann sublevel sets \cite{LibroBallman} for the polydisc
 as ``generalized horospheres" in the Julia's Lemma.
The Busemann sublevel sets are also used to define the analogous of the Koranyi regions.
In contrast with what happens for the existing generalizations of the JWC's Theorem, in our statement
 the class of admissible directions
at a point, $p,$ of the boundary of the polydisc depends upon an entire family of complex geodesics
\cite{LibroAbate} ``passing through"
$p.$

In order to avoid technical complications and give a more geometric
approach, we deal only with the bidisc $\Delta^2$
and use the following terminology.
The symbol
$K_{\Delta^2}$ will denote the Kobayashi distance on $\Delta^2$.
A map  $\Psi\in
\Hol(\Delta,\Delta^{2})$ is called a complex geodesic passing through $y\in\overline{\Delta^2}$
if it is an isometry between the Poincar\'e distance
$\omega$ of $\Delta$ and the Kobayashi distance of $\Delta^2$ whose
image closure contains the point $y$. Let $\Psi$ be a complex geodesic passing through a
point $y=(y_{1}, y_{2})\in\partial\Delta^2$ and suppose $|y_1|=1$.
 It is well known that at least one
component of $\Psi$ is an automorphism of $\Delta$ (see
\cite[Proposition 2.6.10]{LibroAbate}). Thus, up to
re-parametrization, we can assume $\Psi$ is given by $ \Delta\ni z
\to(z, g(z)), $  for some $ g\in \Hol(\Delta, \Delta)$. Thought it
is not necessary a priori, we will assume that $g$ has some
(natural) regularity at $y_1$, namely that the non-tangential limit of
$g$ at $y_1$ is $y_2.$ We denote by $\lambda_g(y_{1}):=\lambda_{g}\in (0,+\infty]$
the boundary dilation coefficient of $g$ at $y_1$ (see section $2$).

Since the function $[K_{\Delta^{2}}(x,\Psi(r))- K_{\Delta^{2}}(\Psi(0),\Psi(r))]$
is non increasing and bounded from below,
the \emph{Busemann function}  $B^{\Psi}(x)$  associated to the
geodesic $\Psi$ can be defined by
$$B^{\Psi}(x):=\lim\limits_{r\to 1^-} [K_{\Delta^{2}}(x,\Psi(r))-
K_{\Delta^{2}}(\Psi(0),\Psi(r))]$$ (see e.g. \cite{LibroBallman}
pag.$23$).
The \emph{ Busemann sublevel set of center $y\in\partial\Delta^2$ and radius $R>0$} of the
function $B^{\Psi}(x)$ will be the set
\begin{equation}\label{definiz.sottolivelli di Busemann}
\B^{\Psi}(y,R):=\{x\in\Delta^2:B_{\Psi}(x)\leq \frac{1}{2}\log R\}.
\end{equation}

To begin with  we
study the geometry of the sets $\B^{\Psi}(y,R)$ (see Proposition \ref{busemann sublevels}).
It turns out that for any point $y=(y_{1},y_{2})\in(\partial\Delta)\times(\partial\Delta)$ on the \v{S}ilov
boundary of the bidisc, given any $R>0,$ we have a continuous family of Busemann sublevel sets
of  radius $R$ of the form $$E_{\Delta}(y_{1},R)\times
E_{\Delta}(y_{2},S).$$
Notice  that every product of horocycles can be seen as a Busemann sublevel set in at least two different
ways, indeed:
$$E_{\Delta}(y_{1},R)\times
E_{\Delta}(y_{2},S)=\B^{\varphi_{g}}(y, R)$$ with $\varphi_{g}(z)=(z,g(z))$ and $\lambda_{g}(y_{1})=\frac{S}{R}$
and moreover
$$E_{\Delta}(y_{1},R)\times
E_{\Delta}(y_{2},S)=\B^{\varphi_{h}}(y, S)$$ with $\varphi_{h}(z)=(h(z),z)$ and $\lambda_{h}(y_{2})=\frac{R}{S}.$
From now on we will simply write $$\B(y,R)=E_{\Delta}(y_{1},R)\times E_{\Delta}(y_{2},R)$$ to denote the Busemann sublevel set associated
to a complex geodesic $\varphi_{g}$ passing through $y$ such that $\lambda_{g}=1.$
Instead, for the points on the flat components of the boundary of $\Delta^2$ there is only one
Busemann sublevel set of a given radius $R>0$ of the form
$$\Delta\times
E_{\Delta}(y_{2},S)\:\: \hbox{if}\:\: (y_{1},y_{2})\in\Delta\times\partial\Delta$$
or
$$E_{\Delta}(y_{1},R)\times
\Delta\:\: \hbox{if}\:\: (y_{1},y_{2})\in\partial\Delta\times\Delta.$$
In the sequel, we will denote by
$\B_{(\lambda_1,\lambda_2)}(y,R)$ (with $\lambda_{1},\lambda_{2}>0$, possibly $+\infty$)
 the Busemann sublevel set given by
$E_{\Delta}(y_{1},\lambda_1 R)\times E_{\Delta}(y_{2},\lambda_{2}R),$
with the convention that
 $E_{\Delta}(y_{i},\lambda_{i}R)=\Delta$ if either
$y_{i}\in\Delta$ or $y_{i}\in\de\Delta$ and $\lambda_{i}=+\infty.$

Our first result is the following version of Julia's lemma:

\begin{theorem}\label{JuliaBusemannlemma}
Let $f=(f_{1},f_{2})\in\Hol(\Delta^2,\Delta^2)$. Let $x=(x_{1},
x_{2})\in\partial(\Delta\times\Delta)=\partial\Delta^2$
and let (for example) $\varphi_{g}(z)=(z,g(z))$ be a complex
geodesic passing through $x.$
Let
$$\frac{1}{2}\log\lambda_{j}:=\lim\limits_{t\to
1^-}[K_{\Delta^2}(0,\varphi_{g}(tx_{1}))
-\omega(0,f_{j}(\varphi_{g}(tx_{1}))]\:\:\: j=1,2$$
Suppose that either $\lambda_{1}<\infty$ or $\lambda_{2}<\infty.$ Then
there exists  a point $y=(y_{1},y_{2})\in \partial\Delta^2$  such
that  for all $R>0$
$$f(\B_{(1,\lambda_{g})}(x,R))\subseteq\B_{(\lambda_{1},\lambda_{2})}(y,R).$$
\end{theorem}
A second achievement is the proof of a generalization of the Lindel\"{o}f Theorem which is based
on the definition of admissible limits. Let
$x\in\partial\Delta^2.$ A continuous curve $\sigma(t)\subset
\Delta^2$ converging to $x$ as $t\to 1^-$ is called a $x-$curve. Let
$\varphi_{g}:\Delta\to \Delta^2$ be a complex geodesic passing
through $x$ and parameterized by $z\mapsto (z, g(z))$ with
$g\in\Hol(\Delta,\Delta).$ A holomorphic function
${\tilde{\pi}_{g}}:\Delta^2 \to \Delta$   such that:
${\tilde{\pi}_{g}}\circ\varphi_{g}=\id_{\Delta}$ is called a
\textsl{g-left inverse} of $\varphi_{g}$. The composition
$\pi_{g}:=\varphi_{g}\circ{\tilde{\pi}_{g}}:\Delta^2\to\Delta^2,$
(such that $\pi_{g}\circ\varphi_{g}=\varphi_{g},$ and
${\tilde{\pi}_{g}}\circ{\tilde{\pi}_{g}}={\tilde{\pi}_{g}}$) is
called a \textsl{g-holomorphic retraction}. The pair
$(\varphi_{g},  \pi_{g})$ is a \textsl{g-projection device}.
Existence of $g$-projection devices, also known as Lempert's
projection devices, in convex domains is established in \cite{RoydenWong} (see also \cite{Lempert82}, \cite{LibroAbate}). In strongly convex domains the Lempert's
projection devices are essentially unique (see \cite{BPT}) while in
the bidisc various holomorphic retractions with different ``fibers" may correspond to a given complex geodesic.
The following definitions will be used in the statements of the generalizations of the Lindel\"{o}f Theorem
and JWC's Theorem.

\begin{definition}
Let $x\in\partial\Delta^2$ and $M>1.$ The \emph{g-Koranyi region}
 $H_{\varphi_{g}}(x,M),$ of vertex $x$ and amplitude $M$ is:
 \begin{equation*}
 \begin{split}
 H_{\varphi_{g}}(x,M):=\{z\in\Delta^2 :\: \: \lim\limits_{r\to 1^-} K_{\Delta^{2}}(z,\varphi_{g}(r))-
&K_{\Delta^{2}}(\varphi_{g}(0),\varphi_{g}(r))\\& + K_{\Delta^{2}}(\varphi_{g}(0),z)<\log M\}.
\end{split}
\end{equation*}
A holomorphic function $f\in\Hol(\Delta^2,\Delta)$ has {\sl $K_{g}-$limit}  equal to $L\in\C$ if $f$ approaches to $L$
inside any  g-Koranyi region.
\\The function $f$ is {\sl
$K_{g}-$bounded} if $\forall $ $M$ there exists a constant $C_{M}>0$ such that
$||f(z)||<C_{M}$ for all $z\in H_{\varphi_{g}}(x,M). $
\end{definition}

\begin{definition}\label{nuovaDefinizione specialita' e ristrettezza}
Let $\sigma(t)\subset \Delta^2$ be a $x-$curve.
\begin{itemize}
\item[-] the curve $\sigma(t)$ is {\sl $g$-special} if
$K_{\Delta^2}(\sigma(t), \pi_{g}(\sigma(t)))\to 0$ as $t\to 1^-$.
\item[-] the curve $\sigma(t)$ is {\sl $g$-restricted } if  ${\tilde{\pi}_{g}}(\sigma(t))\to
\tilde{\pi}_{g}(x)$  non-tangentially as $t\to 1^-$, for $j=1,2.$
\end{itemize}
Moreover, if $h:\Delta^n \to \C$ is holomorphic we say that $h$ has
{\sl restricted $K_{g}-$limit} equal to $L\in\C$ if $h$ has limit
$L$ along any curve which is {\sl$g$-special} and {\sl
$g$-restricted,} and we write
\[
\displaystyle\tilde{K_{g}}-\!\lim_{\hspace{-.5cm}z\to x}h(z)=L .\]
\end{definition}
The announced Lindel\"{o}f type Theorem, proved in this paper, has the following statement:
\begin{theorem}\label{LindeloffBusemannTheorem}
Let $f\in\Hol(\Delta^2,\Delta)$ be a holomorphic function. Given
$x\in\partial\Delta^2$ let $\varphi_{g}$ be a complex geodesic
passing through $x.$ Assume that  $f$ is
$K_{g}-$bounded. If $\sigma_{0}$ is a  $g$-special and $g$-restricted
 $x-$curve such that
$$\lim\limits_{t\to 1^-}f(\sigma_{0}(t))=L$$
then $f$ admits restricted $K_{g}-$limit equal to $L$ at $x.$
\end{theorem}
The above result plays a key role in the proof of our main result:
\begin{theorem}\label{JWCBusemannTheorem}
Let $f\in\Hol(\Delta^2,\Delta^2)$  and $x\in\partial\Delta^2.$
Let  $\varphi_{g}$ be any complex geodesic  passing through $x$ and
parameterized by $\varphi_{g}(z)=(z,g(z)),$ with
$g\in\Hol(\Delta,\Delta).$
Let
$\tilde{\pi_{g}}:\Delta^2\to\Delta$ be the
 g-left inverse
of $\varphi_{g}$ given by $\tilde{\pi_{g}}(z_{1},z_{2})=z_{1}.$
Suppose that for $j=1,2$
$$\frac{1}{2}\log\lambda_{j}=\lim\limits_{t\to 1}[K_{\Delta^2}(0,\varphi_{g}(tx_{1}))-\omega(0,f_{j}
(\varphi_{g}(tx_{1})))]<\infty.$$ Then
there exists a point $y=(y_{1},y_{2})\in(\partial\Delta)^2$ such that the restricted
$K_{g}-$limit of $f_{j}$ at $x$ is $y_{j}$
  for $j=1,2,$ and
$$\displaystyle\tilde{K_{g}}-\!\lim_{\hspace{-.5cm}z\to x}\frac{y_{j}-f_{j}(z)}{1-\tilde{\pi_{g}}(z)}
=\lambda_{j}\min\{1,\lambda_{g}\}$$
$$\displaystyle\tilde{K_{g}}-\!\lim_{\hspace{-.5cm}z\to x}\frac{y_{j}-f_{j}(z)}{1-z_{2}}=\frac{\lambda_{j}}{\max\{1,\lambda_{g}\}}.$$
\end{theorem}

The paper is organized as follows:  in Section $2$ we study in
detail the geometry  the Busemann sublevel sets. In Section $3$
we discuss of special and restricted curves. In Section $4$ we
introduce a new extension of the notion of non-tangential limits and in Section
$5$ we prove  a new version of the Lindel\"of Theorem. In Section $6$ we give
 our extension  of the classical Julia's
Lemma. In Section $7$ we prove our generalization of the
Julia-Wolff-Carath\'{e}odory Theorem. We end the paper in Section
$8$ with an application of our results to the study of the dynamics of fixed points
free holomorphic self-maps of the bidisc.
In fact in this section we give a geometrical
interpretation of a result due to Herv\'{e} \cite{HerveSuBidisco} in terms of the set of generalized Wolff points of
a fixed point free $f\in\Hol(\Delta^2,\Delta^2).$

I would like to sincerely thank Stefano Trapani for his useful suggestions.  Filippo Bracci for his support
and many helpful conversations and Graziano Gentili whose many comments improved this work.

\section{Busemann Functions and a family of horospheres}
The aim of this section is to study in detail the Busemann
sublevel sets and their relation with the horospheres in the
polydisc. Let $\varphi_{g}\in \Hol(\Delta,\Delta^{2})$ be a
complex geodesic, passing through a point $y=(y_{1},
y_{2})\in\partial\Delta^2,$ parameterized as $z\to(z,g(z)).$
Let denote by $\lambda_{g}$ the \emph{boundary dilation coefficient} of $g$ at $y_{1},$
that is
\[\lambda_{g}:= \liminf_{z\to
 y_{1}}\frac{1-|g(z)|}{1-|z|}.\]

 In strongly convex
domains of $\C^n$ the definition of Busemann sublevel set is equivalent to the definition
of horosphere, but in the bidisc this is no longer true. Let
$E(y,R)$ be the \emph{small horosphere}   of center $y\in\Delta^n$
and radius $R$ given by
$$  E(y,R)=\left\{z\in D:\limsup\limits_{w\to y}
  [K_{\Delta^n}(z,w)-K_{\Delta^n}(0,w)]<\frac{1}{2}\log R\right\},$$
and let  $F(y,R)$ be the \emph{big horosphere}   of center $y$ and radius $R$ given by
$$F(y,R)=\left\{z\in D:\liminf\limits_{w\to y} [K_{\Delta^n}(z,w)-K_{\Delta^n}(0,w)]<\frac{1}{2}\log R\right\}.$$
If $n=1$ then $F(y,R)\equiv E(y,R)\equiv E_{\Delta}(y,R)\subset \Delta,$  the horocycle centered in $y$ with radius $R.$
If $n>1 $ then the following proposition holds:
\begin{proposition}\label{busemann sublevels}
Let $\varphi_{h}(z)=(\theta(z), h(z))$ be a complex geodesic in
$\Delta^2$ passing through a point
$y\in\partial\Delta^2,$ where $\theta\in\Aut(\Delta)$ and $h\in\Hol(\Delta,\Delta).$

If $y=(e^{\imath\alpha_{1}},e^{\imath\alpha_{2}})\in(\partial\Delta)^2,$ let  $\lambda_{\theta}$ and
$\lambda_{h}$ respectively be the boundary dilation coefficients of the maps $\theta$ and $h,$ at $e^{\imath\alpha_{1}}$
then
$$\B^{\varphi_{h}}(y,R)=E_{\Delta}(e^{\imath\alpha_{1}},\lambda_{\theta}R)\times
E_{\Delta}(e^{\imath\alpha_{2}},\lambda_h R)$$
If $y\in[\partial\Delta^2 \setminus (\partial\Delta)^2]$
then
$$\B^{\varphi_{g}}(y,R)\equiv E(y,R)\equiv F(y,R).$$
\end{proposition}
\begin{proof}
Up to conjugation with automorphisms,  we can suppose
$(e^{\imath\alpha_{1}},e^{\imath\alpha_{2}})=(1,1).$ Then
$h(1)=1,$ in the sense of non-tangential limit and $\theta(1)=1$.
Let first suppose that $x\in\B^{\varphi_{h}}((1,1),R),$ then, by definition of Busemann sublevel sets,
$\lim\limits_{r\to 1}
[\max\{\omega(x_{1},\theta(r)), \omega(x_{2}, h(r))\}-\omega(0,r)]\leq
\frac{1}{2}\log R.$
We consider the two following cases:

$a)$ there exists a sequence $\{r_{k}\}_{k\in\N}\subseteq (0,1)$ such that
$r_{k}\to 1^-$ as $k\to\infty$ and
$\max\{\omega(x_{1},\theta(r_{k})),\omega(x_{2},
h(r_{k}))\}=\omega(x_{1},\theta(r_{k})),$ and

$b)$ there exists a sequence $\{r_{k}\}_{k\in\N}\subseteq (0,1)$ such that
$r_{k}\to 1^-$ as $k\to\infty$ and
$\max\{\omega(x_{1},\theta(r_{k})),\omega(x_{2},
h(r_{k}))\}=\omega(x_{2},h(r_{k})).$

In case $a)$ we have that
\begin{equation*}
\begin{split}
\frac{1}{2}&\log R \geq \lim\limits_{r\to 1}[
\max\{\omega(x_{1},\theta(r)), \omega(x_{2},
h(r))\}-\omega(0,r)]\\
&=\lim\limits_{k\to\infty} [\omega(x_{1}, \theta(r_{k}))-\omega(0,r_{k})]=\\
&=\lim\limits_{k\to\infty} [\omega(x_{1}, \theta(r_{k}))-\omega(0, \theta(r_{k}))
+\omega(0, \theta(r_{k}))-\omega(0,r_{k})]=\\
&=\lim\limits_{r\to 1^-} [\omega(x_{1}, \theta(r))-\omega(0, \theta(r))
+\omega(0, \theta(r))-\omega(0,r)]=\\
&=\lim\limits_{w\to 1^-} [\omega(x_{1}, w)-\omega(0, w)]+\frac{1}{2}\log \frac{1}{\lambda_{\theta}}.
\end{split}
\end{equation*}
It follows
that $x_{1}\in E_{\Delta}(1,\lambda_{\theta} R).$
Moreover
\begin{equation*}
\begin{split}
&\frac{1}{2}\log R \geq\lim\limits_{k\to\infty} [\omega(x_{1}, \theta(r_{k}))-\omega(0,r_{k})]\\
&\geq
\lim\limits_{k\to \infty}
[\omega(x_{2},h(r_{k}))-\omega(0,r_{k})]\\
&= \lim\limits_{k\to
\infty} [\omega(x_{2},h(r_k))-\omega(0,h(r_k))]+\lim\limits_{k\to
\infty} [\omega(0,h(r_k))-\omega(0,r_k)]\\
&=\lim\limits_{r\to
1}[\omega(x_{2},h(r))-\omega(0,h(r))]+\lim\limits_{r\to 1} [\omega(0,
h(r))-\omega(0,r)]\\
&=\lim\limits_{w\to
1}[\omega(x_{2},w)-\omega(0,w)]+\frac{1}{2}\log\frac{1}{\lambda_{h}}
\end{split}
\end{equation*}
then $x_{2}\in E_{\Delta}(1,\lambda_{h}R).$
Thus, in case $a),$ the first inclusion is proved and
  $\B^{\varphi_{h}}((1,1),R)\subseteq E_{\Delta}(1,\lambda_{\theta}R)\times E_{\Delta}(1,\lambda_{h}R). $

In case $b)$ we notice that
\begin{equation*}
\begin{split}
&\frac{1}{2}\log R \geq \lim\limits_{r\to 1}
[\max\{\omega(x_{1},\theta(r)), \omega(x_{2},
h(r))\}-\omega(0,r)]
\\&=\lim\limits_{k\to\infty} [\omega(x_{2},
h(r_{k}))-\omega(0,r_{k})]= \lim\limits_{r\to 1^-} [\omega(x_{2},
h(r))-\omega(0,r)]
\\ &
\geq  \lim\limits_{k\to \infty}
[\omega(x_{1},\theta(r_{k}))-\omega(0,r_{k})].
\end{split}
\end{equation*}
then proceeding as in case $a)$ it follows that $x_{1}\in E_{\Delta}(1,\lambda_{\theta}R)$ and
$x_{2}\in E_{\Delta}(1,\lambda_{h}R).$
We conclude that, also in this case, $\B^{\varphi_{h}}((1,1),R)\subseteq E_{\Delta}(1,\lambda_{\theta}R)\times
E_{\Delta}(1,\lambda_{h}R)$.

On the other hand, if a point $x\in E_{\Delta}(1,\lambda_{\theta}R)\times E_{\Delta}(1,\lambda_{h}R)$ then
by definition of horocycle
$$\lim_{w\to 1} [\omega(x_{1}, w)-\omega(0,w)]=\lim_{r\to 1^-}
[\omega(x_{1}, \theta(r))-\omega(0,\theta(r))]\leq \frac{1}{2}\log \lambda_{\theta} R$$
and
$$\lim_{r\to 1} [\omega(x_{2}, h(r))-\omega(0,h(r))]\leq
\frac{1}{2}\log \lambda_{h}R.$$ Thus
\[\lim_{r\to 1^-} [\omega(x_{2}, h(r))-\omega(0,r)]
=\lim_{r\to 1^-}
[\omega(x_{2},
h(r))-\omega(0,h(r))+\omega(0,h(r))-\omega(0,r)]\]

\[\leq
\frac{1}{2}\log \lambda_{h}R\frac{1}{\lambda_{h}} =\frac{1}{2}\log R.\]
Swapping $h$ with $\theta,$ arguing as above, we have
\[\lim_{r\to 1^-} [\omega(x_{1}, \theta(r))-\omega(0,r)]\leq\frac{1}{2}\log R.\]
We conclude that
\begin{equation}\label{limitBusemann}
\lim\limits_{r\to 1} [\max\{\omega(x_{1},\theta(r)), \omega(x_{2},
h(r))\}-\omega(0,r)]\leq \frac{1}{2}\log R
\end{equation}
and $\B^{\varphi_{h}}((1,1),R)= E_{\Delta}(1,\lambda_{\theta}R)\times E_{\Delta}(1,\lambda_{h}R).$

If $y$ is a point of the flat component of $\partial\Delta^2,$ then
$E(y,R)\equiv F(y,R)$ and therefore the limit that defines the
small and big horosphere exists. Thus it follows immediately that $\B^{\varphi}(y,R)\equiv E(y,R)\equiv F(y,R)$
for all geodesic $\varphi$ and for all $R>0.$
\end{proof}

Let us notice that if we consider the re-parametrization of $\varphi_{h}(z)=(\theta(z),h(z))$
given by $\varphi_{g}(z)=(z, g(z))$ where $g:=h\circ\theta^{-1}$ we have
$$\B^{\varphi_{g}}((1,1),R)=E_{\Delta}(1,R)\times
E_{\Delta}(1, \lambda_{g} R),$$ where $\lambda_{g}=\frac{\lambda_{h}}{\lambda_{\theta}}.$
It follows by the same arguments used in Proposition \ref{busemann sublevels}.
For this reason, from now on, we  consider only parametrization of the type $z\to(z,g(z))$ (respectively
$z\to(g(z),z))$ (see also the Introduction).

For later use we now compute explicitly the Busemann sublevel sets. We use the above notation.
Let us first consider a point  $y=(y_{1},y_{2})$  contained in a flat component
of the boundary of the bidisc, $\partial\Delta^2.$
As Proposition \ref{busemann sublevels} states, the Busemann sublevel sets, centered in $y,$ coincide with
the small and the big horosphere (see Abate \cite{ArticoloJWCAbate} for an explicit description os small and big
horospheres).
On the other hand let us consider a point $y=(y_{1},y_{2})\in(\partial\Delta)^2.$
Without loss of generality we can suppose that $y=(1,1).$ We claim that
\begin{equation}\label{forma euclidea sottolivelli}
E_{\Delta}(1,R)\times E_{\Delta}(1,\lambda_{g}R)=\left\{z\in\Delta^2:\max_{j=1,2}
\frac{|1-z_{j}|^2}{(1-|z_{j}|^2)}\lim\limits_{r\to
1}\frac{(1-r^{2})}{ (1-|\varphi_{g_{j}}(r)|^2)}\leq R\right\}.
\end{equation}

Indeed assume  that

\[\max_{j=1,2}
\left\{\frac{|1-z_{1}|^2}{(1-|z_{1}|^2)},
\frac{|1-z_{2}|^2}{(1-|z_{2}|^2)}\frac{1}{\lambda_{g}}\right\}\leq R.\]
Thus the two  possibilities hold:
\\
\\
$i)\:\:\max_{j=1,2} \left\{\frac{|1-z_{1}|^2}{(1-|z_{1}|^2)},
\frac{|1-z_{2}|^2}{(1-|z_{2}|^2)}\frac{1}{\lambda_{g}}\right\}=
\frac{|1-z_{1}|^2}{(1-|z_{1}|^2)}\:\:$ then
$\:\frac{|1-z_{2}|^2}{(1-|z_{2}|^2)}\frac{1}{\lambda_{g}}\leq
\frac{|1-z_{1}|^2}{(1-|z_{1}|^2)}\leq R $
\\
\\
and, by definition of horocycles,  $ z_{1}\in E_{\Delta}(1,R)$
and  $z_{2}\in E_{\Delta}(1,\lambda_{g}R).$
\\
\\
$ii) \:\: \max_{j=1,2} \left\{\frac{|1-z_{1}|^2}{(1-|z_{1}|^2)},
\frac{|1-z_{2}|^2}{(1-|z_{2}|^2)}\frac{1}{\lambda_{g}}\right\}=
\frac{|1-z_{2}|^2}{(1-|z_{2}|^2)}\frac{1}{\lambda_{g}}\:\:$ then
$\frac{|1-z_{1}|^2}{(1-|z_{1}|^2)}\leq\frac{|1-z_{2}|^2}
{(1-|z_{2}|^2)}\frac{1}{\lambda_{g}}\leq \nolinebreak[4]R $
\\
\\
and, by definition of horocycles again  $z_{1}\in
E_{\Delta}(1,R)$ and $z_{2}\in E_{\Delta}(1,\lambda_{g}R).$
Namely $z\in E_{\Delta}(1,R)\times
E_{\Delta}(1,\lambda_{g}R).$ Conversely let
$z=(z_{1},z_{2})\in E_{\Delta}(1,R)\times E_{\Delta}(1,\lambda_{g}R)=\B^{\varphi_{g}}((1,1),R),$ with $\varphi_{g}(z)=
(z,g(z)).$  By definition of
horocycles, it follows that either
$\frac{|1-z_{1}|^2}{(1-|z_{1}|^2)}\leq R$ or $\frac{|1-z_{2}|^2}
{(1-|z_{2}|^2)}\frac{1}{\lambda_{g}}\leq R,$ hence
$$\max_{j=1,2}
\left\{\frac{|1-z_{1}|^2}{(1-|z_{1}|^2)},
\frac{|1-z_{2}|^2}{(1-|z_{2}|^2)}\frac{1}{\lambda_{g}}\right\}\leq R.$$
By the very  definition of sublevel sets of Busemann functions,
$z=(z_{1},z_{2})\in\B^{\varphi_{g}}((1,1),R),$ proving the claim.

\section{Special and restricted curves}

Let be $x=(x_{1},x_{2})\in\partial\Delta^2$ and  $\varphi_{x}:\Delta\to\Delta^2$ the complex geodesic passing through $x,$
defined by
$\varphi_{x}(z)=zx.$
Let us denote by $d_{x}$  the
{\emph{$\check{S}$ilov degree}} of $x,$ that is the number of
components of $x$ with absolute value $1$ and
$\check{x}:=(\check{x_{1}},\check{{x}_{2}})$ is the
$\check{S}$ilov part of $x,$ defined by
$$\check{x}_{j}=\left\{
\begin{array}{lcccc}
x_{j} & \hbox{if}& |x_{j}| & = & 1, \\
0 & \hbox{if}& |x_{j}| & < & 1.
\end{array}\right.
$$
In this setting Abate \cite{ArticoloJWCAbate} gave the following definition
\begin{definition}\label{AbateProjection}

We call the holomorphic function  $\tilde{p}_{x}:\Delta^2\to\Delta$  given by
$$\tilde{p}_{x}(z):=\frac{1}{d_{x}}(z,\check{x}),$$  such that
$\tilde{p}_{x}\circ \varphi_{x}=id_{\Delta^n},$  an \emph{Abate's left
inverse} of $\varphi_{x}.$

We call the holomorphic function  $p_{x}:\Delta^n\to\Delta^n$ given by
 $p_{x}(z):=\varphi_{x}\circ \tilde{p}_{x},$ such that
 $p_{x}\circ p_{x}=p_{x}$ and
 $p_{x}\circ \varphi_{x}=\varphi_{x}$ an \emph{Abate's holomorphic retraction}.

A $x-$curve $\sigma(t)\subset \Delta^n$
 is {\sl A-special} if $
K_{\Delta^n}(\sigma(t), p_{x}(\sigma(t)))\to 0, $ as $t\to 1^-$.

A $x-$curve
$\sigma(t)\subset \Delta^n$
is {\sl A-restricted} if $\tilde{p_{x}}(\sigma(t))$ converges to $1$
non-tangentially.

The pair $(p_{x},\varphi_{x})$ is called an $A-$projection device.
\end{definition}


We notice that the $A-$projection device due to Abate is not unique,
that is, given the complex geodesic $\varphi_{x}=zx,$ the left inverse $\tilde{p_{x}},$ and the
holomorphic retraction $p_{x}$ are not unique.
Moreover $\varphi_{x}$ is not the unique complex geodesic passing through $x.$
Thus we are led to give the following definitions:

\begin{definition}\label{new projection device}
Let $\varphi_{g}:\Delta\to \Delta^2$
be a complex geodesic passing through $x$ and parameterized by $z\mapsto
(z, g(z));\: g\in\Hol(\Delta,\Delta).$

A holomorphic function  ${\tilde{\pi}_{g}}:\Delta^2 \to \Delta$   such that
${\tilde{\pi}_{g}}\circ\varphi_{g}=\id_{\Delta}$
is called a  \textsl{g-left inverse} function of $\varphi_{g}$ .

A holomorphic function  $\pi_{g}:=\varphi_{g}\circ{\tilde{\pi}_{g}}:\Delta^2\to\Delta^2,$
such that $\pi_{g}\circ\varphi_{g}=\varphi_{g},$ and
${\tilde{\pi}_{g}}\circ{\tilde{\pi}_{g}}={\tilde{\pi}_{g}}$ is called a \textsl{g-holomorphic retraction}.

The pair
$(\varphi_{g}, \pi_{g})$ is a \textsl{g-projection
device}.
\end{definition}

\begin{definition}\label{nuovaDefinizione specialita' e ristrettezza}
Let $\sigma(t)\subset \Delta^2$ be a $x-$curve.
\begin{itemize}
\item[-] the curve $\sigma(t)$ is {\sl $g$-special} if
$K_{\Delta^2}(\sigma(t), \pi_{g}(\sigma(t)))\to 0$ as $t\to
1^-$. \item[-] the curve $\sigma(t)$ is {\sl $g$-restricted }
if  ${\tilde{\pi}_{g}}(\sigma(t))\to {\tilde{\pi}_{g}}(x)$  non-tangentially as $t\to 1^-$, for $j=1,2.$
\end{itemize}
\end{definition}

As a matter of notation, when we refer to the geodesic
parameterized by $\varphi(z)= (z, z),$ we omit the index $g,$
since $g=\id_{\Delta}.$ In addition we denote by $p$ the
$A-$holomorphic retraction given in definition
\ref{AbateProjection}.

 In this setting,
an $x-$curve $\gamma$ is A-special and A-restricted
if $K_{\Delta^2}(\gamma(t), p(\gamma(t)))\to 0$ as $t\to
1^-$ and ${\tilde{p}}(\gamma(t))$ approach to the point $\tilde{p}(x)$ non tangentially.

We notice that if $x$ is a point on a flat component of $\partial\Delta^2,$ definitions \ref{AbateProjection} and
\ref{nuovaDefinizione specialita' e ristrettezza} are equivalent, and we have that
a $x-$curve $\gamma$ is $g$-special and $g$-restricted
 if and only if it is A-special and A-restricted.
On the other hand  if $x\in(\partial\Delta)^2,$ we have the following characterization:

\begin{proposition}\label{special1}
Let denote by $\varphi_{x}(z)$ the complex
 geodesic passing through the point $x\in(\partial\Delta)^2$ parameterized by $z\to(zx)$  and
let $\pi_{g}:\Delta^2\to\varphi_{x}(\Delta)$  be any linear holomorphic retraction on the
 image of the complex geodesic $\varphi_{x}.$
 Let $\gamma(t)=(\gamma_{1}(t),\gamma_{2}(t))$ be a $g-$restricted, $x-$curve in $\Delta^2.$
 Then $\gamma$ is $g-$special if and only if $\frac{1-\gamma_{2}(t)}{1-\gamma_{1}(t)}\to 1$
 as $t\to 1^-.$
\end{proposition}
\begin{proof}
Without loss of generality we may suppose $x=(1,1);$ then $\varphi_{x}(z)=\varphi(z)=(z,z).$
Consider a linear projection $\pi_{g}(z_{1}, z_{2})=(a z_{1}+bz_{2}, a z_{1}+bz_{2})$ with $a,b\in\C.$
By definition of holomorphic retraction we have that $\pi_{g}(z, z)=(az+bz, a z+bz)=(z,z),$ then $b=1-a.$
We know that $\gamma $ is $g-$special if and only if $K_{\Delta^2}(\gamma(t), \pi_{g}(\gamma(t)))\to 0$ as
$t\to 1^-.$
By the very definition of the Kobayashi distance, this is equivalent to
\begin{equation*}
\begin{split}
\lim\limits_{t\to 1^-}
\left|\frac{\gamma_{1}(t)-a\gamma_{1}(t)-(1-a)\gamma_{2}(t)}{1-\overline{\gamma_{1}(t)}[a\gamma_{1}(t)+(1-a)\gamma_{2}(t)]}\right|=0
\end{split}
\end{equation*}
By an easy calculation we get:
\begin{equation*}
\begin{split}
\left|\frac{\gamma_{1}(t)-a\gamma_{1}(t)-(1-a)\gamma_{2}(t)}
{1-\overline{\gamma_{1}(t)}[a\gamma_{1}(t)+(1-a)\gamma_{2}(t)]}\right|=
\left|\frac{(1-a)(\frac{1-\gamma_{2}(t)}{1-\gamma_{1}(t)}-1)}
{\frac{1-\overline{\gamma_{1}(t)}[a\gamma_{1}(t)+(1-a)\gamma_{2}(t)]}{1-\gamma_{1}(t)}}     \right|
\end{split}
\end{equation*}
Thus if the curve $\gamma$ is $g-$special we necessarily have that
$(1-\gamma_{2}(t))/(1-\gamma_{1}(t))\to 1$ as $t\to 1^-.$
On the other hand if
\begin{equation}\label{caratter.curve speciali}
\lim\limits_{t\to 1^-}\frac{1-\gamma_{2}(t)}{1-\gamma_{1}(t)}=1,\end{equation}
 taking into account that $\gamma$ is $g-$restricted, then
\begin{equation*}
\begin{split}
&\left|\frac{1-\overline{\gamma_{1}(t)}[a\gamma_{1}(t)+(1-a)\gamma_{2}(t)]}{1-\gamma_{1}(t)}\right|\geq
\frac{1-|a\gamma_{1}(t)+(1-a)\gamma_{2}(t)|}{|1-\gamma_{1}(t)|}\\&\geq
\frac{|1-a\gamma_{1}(t)-(1-a)\gamma_{2}(t)|}{M|1-\gamma_{1}(t)|}=\frac{|a\frac{1-\gamma_{1}(t)}{1-\gamma_{1}(t)}
+(1-a)\frac{1-\gamma_{2}(t)}{1-\gamma_{1}(t)}|}{M}\to \frac{1}{M}
\end{split}
\end{equation*}
as $t\to 1^- ,$  and condition (\ref{caratter.curve speciali}) is also sufficient.
\end{proof}
It is worth noticing that the Abate projection $p$ is a special linear projection with $a=b=\frac{1}{2}:$

\begin{proposition}\label{specialK}
Let $x\in(\partial \Delta)^2$ and let $\gamma(t)=(\gamma_1(t),\gamma_2(t))$ be a $x-$curve in
$\Delta^2.$ Let
$\varphi_{x}(z)$  be the complex geodesic, parameterized by $z\to(zx),$
passing through the point $x.$  Let
$\pi_{g}$ be any linear projection on $\varphi_{x}$.
Then $\gamma$ is $g$-special and $g$-restricted if and only if
$\gamma$ is $A-$special and $A-$restricted.

\end{proposition}

\begin{proof}
Without loss of generality we suppose $x=(1,1)$ and thus $\varphi_{x}(z)=\varphi(z)=(z,z).$
As in the proof of Proposition \ref{special1} we can write that
$\pi_{g}(z_{1},z_{2})=(az_{1}+(1-a)z_{2}, az_{1}+(1-a)z_{2}),$ with $a,b\in\C.$

We first prove the ``$\Leftarrow$" implication.
We show  that $\gamma$ is $g$-special, and in particular  that
$K_{\Delta^2}(\gamma(t),\pi_{g}(\gamma(t)))\to 0$ as $t\to 1^-.$
By the triangular inequality
\begin{equation}\label{disuguaglianzaPROIEZIONI}
\begin{split}
K_{\Delta^2}(\gamma(t),\pi_{g}(\gamma(t)))\leq
K_{\Delta^2}(\gamma(t),p(\gamma(t)))+ K_{\Delta^2}(\pi_{g}(\gamma(t)),p(\gamma(t)))
\end{split}
\end{equation}
and, since $\gamma $ is a $A-$special curve,  $K_{\Delta^2}(\gamma(t),p(\gamma(t)))\to 0$ as $t\to 1^-.$
Moreover
\begin{equation*}
\begin{split}
K_{\Delta^2}(\pi_{g}(\gamma(t)),p(\gamma(t)))=K_{\Delta^2}(\varphi(\tilde{\pi_{g}}(\gamma(t))),\varphi(\tilde{p}(\gamma(t))))
=\omega(\tilde{\pi_{g}}(\gamma(t)), \tilde{p}(\gamma(t))).
\end{split}
\end{equation*}
We claim that
$\omega(\tilde{\pi_{g}}(\gamma(t)), \tilde{p}(\gamma(t)))\to 0$ as $t\to 1^-.$ By the very definition of
Poincar\'{e} metric, this is equivalent to
$$\lim\limits_{t\to 1^-}\frac{|\tilde{p}(\gamma(t))-\tilde{\pi_{g}}(\gamma(t))|}
{|1-\overline{\tilde{p}(\gamma(t))}\tilde{\pi_{g}}(\gamma(t))|}=0.$$
First we notice:
\begin{equation*}
\begin{split}
&
\frac{|\tilde{p}(\gamma(t))-\tilde{\pi_{g}}(\gamma(t))|}
{|1-\overline{\tilde{p}(\gamma(t))}\tilde{\pi_{g}}(\gamma(t))|}
=\frac{|\tilde{p}(\gamma(t))-1+1-\tilde{\pi_{g}}(\gamma(t))|}
{|1-\overline{\tilde{p}(\gamma(t))}\tilde{\pi_{g}}(\gamma(t))|}
\\&\hspace{1,5cm}=\Large{\left|\frac{\frac{1-\tilde{\pi_{g}}(\gamma(t))}{1-\tilde{p}(\gamma(t))}-1}
{\frac{1-\overline{\tilde{p}(\gamma(t))}\tilde{\pi_{g}}(\gamma(t))}{1-\tilde{p}(\gamma(t))}}\right|}
.\end{split}
\end{equation*}
Moreover, by definition of $A-$projection device,
\[\frac{1-\tilde{\pi_{g}}(\gamma(t))}{1-\tilde{p}(\gamma(t))}=\frac{1-a \gamma_{1}-(1-a)\gamma_{2}}{1-
\frac{1}{2}(\gamma_{1}+\gamma_{2})}\]\[=2\frac{1-a \gamma_{1}-(1-a)\gamma_{2}}{2-\gamma_{1}-\gamma_{2}}=
2\frac{a+(1-a)-a \gamma_{1}-(1-a)\gamma_{2}}{1-\gamma_{1}+1-\gamma_{2}}\]\[
=2\frac{a(1-\gamma_{1})+(1-a)(1-\gamma_{2})}{1-\gamma_{1}+1-\gamma_{2}}=
2\frac{a+(1-a)\frac{(1-\gamma_{2})}{(1-\gamma_{1})}}{1+\frac{(1-\gamma_{2})}{(1-\gamma_{1})}}\]
and by Proposition \ref{special1} we get
\begin{equation*}
\begin{split}
\lim\limits_{t\to 1^-}2\frac{a+(1-a)\frac{(1-\gamma_{2})}{(1-\gamma_{1})}}{1+\frac{(1-\gamma_{2})}{(1-\gamma_{1})}}=1.
\end{split}
\end{equation*}
Furthermore, since by hypothesis the curve $\gamma$ is $A-$restricted, then there exists $M>1$ such that
$$\frac{|1-\gamma(t)|}{1-|\gamma(t)|}<M$$ and in particular
\begin{equation*}
\begin{split}
\left|\frac{1-\overline{\tilde{p}(\gamma(t))}\tilde{\pi_{g}}(\gamma(t))}{1-\tilde{p}(\gamma(t))}\right|\geq
\frac{1-|\tilde{p}(\gamma(t))| |\tilde{\pi_{g}}(\gamma(t))|}{|1-\tilde{p}(\gamma(t))|}\geq
\frac{1-|\tilde{p}(\gamma(t))|}{|1-\tilde{p}(\gamma(t))|}>\frac{1}{M}.
\end{split}
\end{equation*}
Then we conclude that
\begin{equation}\label{distanzaTRA2proiezioni}
\lim\limits_{t\to 1^-}\omega(\tilde{\pi_{g}}(\gamma(t)), \tilde{p}(\gamma(t)))= 0
\end{equation}
 and the curve $\gamma$
is $g$-special. Let prove that it is also $g$-restricted. Let notice that
by equation (\ref{distanzaTRA2proiezioni}) and since the curve $\gamma$ is $A-$restricted
\begin{equation*}
\begin{split}
\frac{|1-\tilde{\pi_{g}}(\gamma(t))|}{1-|\tilde{\pi_{g}}(\gamma(t))|}=
\frac{|1-\tilde{\pi_{g}}(\gamma(t))|}{|1-\tilde{p}(\gamma(t))|}\:
\frac{|1-\tilde{p}(\gamma(t))|}{1-|\tilde{p}(\gamma(t))|}\:\frac{1-|\tilde{p}(\gamma(t))|}
{1-|\tilde{\pi_{g}}(\gamma(t))|}< 4M
\end{split}
\end{equation*}
and then
the curve $\gamma$ is $g$-restricted.

The last step consists in proving the ``$\Rightarrow$" implication of the theorem.
To do this it is sufficient to interchange the Abate's projection $p$ with the linear projection $\pi_{g}$ in the proof above
and the thesis easily follows.
\end{proof}
\begin{remark} By Proposition \ref{specialK} it follows  the Abate's Julia-Wolff-Carath\'{e}odory theorem
for linear projections.
\end{remark}
The next question is what happens if we consider another geodesic passing through the point $x\in(\partial\Delta)^2.$
Arguing as in Proposition \ref{specialK} we have:
\begin{proposition}\label{caratterizzazioneVarphiSPECIALI}
Let $(\varphi_{g}, \pi_{g})$ be a projection
device. Let assume that the geodesic $\varphi_{g}$ passes through a point $x\in(\partial\Delta)^2$ and set
$\lambda_{g}:=\liminf_{z\to x_{1}}\frac{1-|g(z)|}{1-|z|}<\infty.$
Let $\gamma:=(\gamma_{1},\gamma_{2})$ be a $g$-restricted $x-$curve in $\Delta^2.$
Then $\gamma$ is $g$-special if and only
if $$\lim\limits_{t\to 1^-}\frac{1-\gamma_{2}(t)}{1-\gamma_{1}(t)}=\lambda_{g}$$
\end{proposition}
\section{The non-tangential limit}
The  non-tangential limit in $\Delta$ can be defined in two equivalent ways.
We can say that a function $f\in\Hol(\Delta,\Delta)$ has \emph{non-tangential} limit $L\in\C$ at a point
$y\in\partial\Delta$ if $f(z)\to L$ as $z\to y,$ inside any \emph{Stolz region}, $H(y, M)$ of vertex
$y$ and amplitude $M>1,$ where
$$H(y, M):=\left\{z\in\Delta\::\:\frac{|y-z|}{1-|z|}<M   \right\}.$$
We can equivalently  say that $f\in\Hol(\Delta,\Delta)$ has \emph{non-tangential} limit $L\in\C$ at a point
$y\in\partial\Delta$ if $f(\sigma(t))\to L$ as $t\to 1,$ along any curve $\sigma:[0,1)\to\Delta$ such that
$\sigma(t)\to y$ non-tangentially as $t\to 1^-.$
In \cite{ArticoloJWCAbate} (see also \cite{ArticoloLindelofAbate}) Abate generalizes the Stolz
region giving the following definition of (small)\emph{Koranyi region} (of vertex $y\in\partial\Delta^2$ and amplitude $M$),
$$H(y,M):=\{z\in\Delta^2 :\: \: \limsup\limits_{w\to y} K_{\Delta^{2}}(z,w)-
K_{\Delta^{2}}(0,w) + K_{\Delta^{2}}(0,z)<\log M\}.$$
Thus an extension of the first definition of non-tangential limit becomes (see \cite{ArticoloJWCAbate}):
\begin{definition}
A map $f:\Delta^2\to\C^m$ has $K-$limit $L\in\C^m$ at $y\in\partial\Delta^2$ if $f(z)\to L$ as $z\to y$ inside
any Koranyi region.
\end{definition}
On the other hand, by means of $A-$special and $A-$restricted curves, Abate says that (\cite{ArticoloJWCAbate})
a holomorphic
function $f:\Delta^2 \to \C^m$ has {\sl
restricted $K-$limit} $L$ at $x$ if $f(\sigma(t))\to L$ for any
A-special and A-restricted $x-$curve $\sigma(t)\subset \Delta^2$
and we write $\displaystyle\tilde{K}-\!\lim_{\hspace{-.7cm}z\to x} f(z)=L.$

We notice that the definitions of $K-$limit and restricted $K-$limit   are no more equivalent.
More precisely if $f$ has $K-$limit $L$ at $y\in\partial\Delta^2$ then it has restricted $K-$limit too.
The converse is false (see example in \cite{ArticoloJWCAbate}).
We extend these definitions by means of Busemann functions.
The first step consists in giving the following extension of the notion of Stolz region:
\begin{definition}
Let $x\in\partial\Delta^2$ and $M>1,$ the \emph{g-Koranyi region}
 $H_{\varphi_{g}}(x,M),$ of vertex $x$ and amplitude $M$ is:
 \begin{equation}
 \begin{split}
 H_{\varphi_{g}}(x,M):=\{z\in\Delta^2 :\: \: \lim\limits_{r\to 1^-} K_{\Delta^{2}}(z,\varphi_{g}(r))-
&K_{\Delta^{2}}(\varphi_{g}(0),\varphi_{g}(r))\\& + K_{\Delta^{2}}(\varphi_{g}(0),z)<\log M\}.
\end{split}
\end{equation}
\end{definition}

And then we, naturally, say that
\begin{definition}\label{defKglimit}
A holomorphic function $f\in\Hol(\Delta^2,\Delta^2)$ has {\sl $K_{g}-$limit} $L\in\C$ if $f$ approaches to $L$
inside any  g-Koranyi region.
\end{definition}
If we consider the complex geodesic $\varphi(z)=(z,z)$ then the Koranyi region $H((1,1), M)$
coincide with the $g-$Koranyi region $H_{\varphi}((1,1),M)$.

Moreover let $(\varphi_{g}, \pi_{g})$ be a $g-$projection device as in Definition \ref{new projection device}:
\begin{definition}\label{defvarphiglimit}
A holomorphic function $h:\Delta^2 \to \C$ is said to have  {\sl restricted $K_{g}-$  limit} $L$
if $h$ has limit $L$ along any curve which is $g$-special
and {\sl $g$-restricted,} and we write
\[
\displaystyle\tilde{K_{g}}-\!\lim_{\hspace{-.5cm}z\to x}h(z)=L
.\]
\end{definition}
Obviously Definition \ref{defKglimit}  and Definition \ref{defvarphiglimit} ar not equivalent but
again the $K_{g}-$limit implies the restricted $\tilde{K}_{g}-$limit.

\section{Lindel\"{o}f Theorems}

The classical Lindel\"{o}f principle implies that if $f\in\Hol(\Delta,\Delta)$ has limit $L$
 along any given $1-$curve, then $L$ is the non-tangential limit of $f$ at $1.$
The first step to generalize this theorem to several complex variables
consists in detecting a correct class of curves.
Let $(\varphi_{g}, \pi_{g})$
a $g-$projection device.
The idea is to consider the $g$-special and $g$-restricted curves.

In this setting we prove the
following first generalization of the Lindel\"{o}f principle:
\begin{theorem}\label{Lindelof1}
Let $f\in\Hol(\Delta^2, \C)$ be a bounded holomorphic function.
Let  $x\in\partial\Delta^2.$ Assume there exists a $g$-special $x-$curve $\sigma_{0}$
 such that
$$\lim\limits_{t\to 1^-}f(\sigma_{0}(t))=L\in\C.$$
Then $f$ has restricted  $\tilde{K_{g}}-$limit $L$
at $x.$
\end{theorem}
\begin{proof}
This proof is similar to the one  in \cite{ArticoloJWCAbate} (see
Theorem $2.1$). We first observe that, given $\sigma$ a
$g$-special $x-$curve,
\begin{equation}\label{osservazioneLindelof1}
\begin{split}
0\leq \omega(f(\sigma(t)),f(\pi_{g}(\sigma(t))))\leq K_{\Delta^2}(\sigma(t),\pi_{g}(\sigma(t)))\to 0
\end{split}
\end{equation}
as $t\to 1^-.$ Therefore the limit of $f(\pi_{g}(\sigma(t)))$
exists, as $t\to 1^-,$ if and only if the limit of $f(\sigma(t))$
as $t\to 1^-$ does, and the two limits are equal. In particular
$f(\pi_{g}(\sigma_{0}(t)))\to L$ as $t\to 1^-$ and by classical
Lindel\"{o}f principle $f(\pi_{g}(\sigma(t)))\to L$ for any
$g$-restricted $x-$curve and by remark
(\ref{osservazioneLindelof1}) it follows that $f(\sigma(t))\to L$
for any $g$-restricted and $g$-special $x-$curve $\sigma$.
\end{proof}
Since in
 the Julia-Wolff-Carath\'{e}odory theorem the functions we deal with are incremental ratios,
a stronger result than Theorem \ref{Lindelof1} is needed.
It is worthwhile to introduce some definitions and preliminary results.
\begin{definition}
Let $f\in\Hol(\Delta^2, \C).$ We say that $f$ is {\sl
$K_{g}-$bounded} if $\forall $ $M$ there exists a constant $C_{M}>0$ such that
$||f(z)||<C_{M}$ for all $z\in H_{\varphi_{g}}(x,M). $
\end{definition}
\begin{lemma}\label{lemma1preparazioneLindeloff2}
Let $x\in\partial\Delta^2$ and let $(\varphi_{g},\pi_{g})$ a $g-$projection device.
Suppose $\sigma$  an $x-$curve. Then
$\sigma$ is $g$-restricted if and only if
$\pi(\sigma(t))\in H_{\varphi_{g}}(x,M)$ eventually.

\end{lemma}
\begin{proof}
The proof follows by definition of $g$-restricted curve.
Indeed, since $\varphi_{g}$ is a geodesic and
$\pi_{g}(\sigma(t))=\varphi_{g}\circ\tilde{\pi_{g}}(\sigma(t))$
\begin{equation}
\begin{split}
\lim\limits_{w\to x} \omega(\tilde{\pi_{g}}(\sigma(t)),w)-\omega(0,w)+\omega(0,\tilde{\pi_{g}}(\sigma(t)))<\log M
\end{split}
\end{equation}
if and only if
\begin{equation}
\begin{split}
&\lim\limits_{s\to 1^-} K_{\Delta^2}(\pi_{g}(\sigma(t)),\varphi_{g}(s))-\omega(0,s)+K_{\Delta^2}(\varphi_{g}(0),\pi_{g}(\sigma(t)))
\\& =\lim\limits_{w\to x} K_{\Delta^2}(\varphi_{g}(\tilde{\pi_{g}}(\sigma(t))),\varphi_{g}(s))
-\omega(0,s)+K_{\Delta^2}(\varphi_{g}(0),\varphi_{g}(\tilde{\pi_{g}}(\sigma(t))))
\\& =\lim\limits_{w\to x} \omega(\tilde{\pi_{g}}(\sigma(t)),w)-\omega(0,w)+\omega(0,\tilde{\pi_{g}}(\sigma(t)))<\log M
\end{split}
\end{equation}

\end{proof}
\begin{remark}\label{comeSCRIVOunaCURVAspeciale}
Consider $\sigma$ a $g$-special $x-$curve. We notice that it is possible to
write $\sigma(t):=(\sigma_{1}(t),\sigma_{2}(t))=\pi_{g}(\sigma(t))+\alpha(t)$ with $\alpha(t):=
(\alpha_{1}(t), \alpha_{2}(t))\to (0,0)$ as $t\to 1^-.$
By definition of the projection $\pi_{g}$ we get that $\alpha_{1}(t)\equiv 0$ and $\alpha_{2}(t)\to 0,$
as $t\to 1^-.$
\end{remark}

\begin{lemma}\label{lemma2preparazioneLindeoff2}
Let $x\in\partial\Delta^2$ and $(\varphi_{g},\pi_{g})$ a $g-$projection device. Let $\sigma$ be an $x-$curve. Write
$\sigma(t)=\pi_{g}(\sigma(t))+\alpha(t)$ with $\alpha(t)\to 0$ as
$t\to 1^-.$ Then $\sigma$ is $g$-special if
and only if
$$\lim\limits_{t\to
1^-}\frac{|\alpha_{2}(t)|}{1-|g(\sigma_{1}(t))|}=0.$$
\end{lemma}
\begin{proof}
Assume  first that
$$\lim\limits_{t\to 1^-}\frac{|\alpha_{2}(t)|}{1-|g(\sigma_{1}(t))|}=0.$$
By the triangular inequality and by definition of Kobayashi distance in the bidisc, we have that
\begin{equation*}
\begin{split}
&K_{\Delta^2}(\sigma(t),\pi_{g}(\sigma(t)))=\max\{\omega(\sigma_{1}(t),\sigma_{1}(t));\omega(\sigma_{2}(t),g(\sigma_{1}(t)))\}\\
& =\frac{1}{2}\log \frac{1+\left|\frac{\sigma_{2}(t)-g(\sigma_{1}(t))}{1-\overline{\sigma_{2}(t)}
g(\sigma_{1}(t))}\right|}{1-\left|\frac{\sigma_{2}(t)-g(\sigma_{1}(t))}{1-\overline{\sigma_{2}(t)}
g(\sigma_{1}(t))}\right|}
\leq \frac{1}{2}\log \frac{1+\frac{|\alpha_{2}(t)|}{1-|g(\sigma_{1}(t))|}}
{1-\frac{|\alpha_{2}(t)|}{1-|g(\sigma_{1}(t))|}}\to 0
\end{split}
\end{equation*}
as $t\to 1^-$. Thus the curve $\sigma$ is $g$-special and the first implication has been proved.
On the other hand, let suppose that $\sigma$ is $g$-special.

If, by contradiction, $\lim_{t\to
1^-}\frac{|\alpha_{2}(t)|}{1-|g(\sigma_{1}(t))|}\neq 0$  then
there  exists $\tilde{\varepsilon}>0$ such that
$\frac{|\alpha_{2}(t)|}{1-|g(\sigma_{1}(t))|}>\tilde{\varepsilon}>0.$
In particular there  exists $\varepsilon>0$ such that
\[T:=\frac{|\alpha_{2}(t)|}{1-|g(\sigma_{1}(t))|^2}>\varepsilon>0.\]
Furthermore
\begin{equation*}
\begin{split}
&\frac{|1-\overline{g(\sigma_{1}(t))}\sigma_{2}(t)|}{1-|g(\sigma_{1}(t))|^2}
=\frac{|1-\overline{g(\sigma_{1}(t))}(g(\sigma_{1}(t))+\alpha_{2}(t))|}{1-|g(\sigma_{1}(t))|^2}=
\left|1-\frac{\overline{g(\sigma_{1}(t))}\alpha_{2}(t)}{1-|g(\sigma_{1}(t))|^2}\right|\\
&\hspace{1cm}\leq 1+\left|\frac{\overline{g(\sigma_{1}(t))}\alpha_{2}(t)}{1-|g(\sigma_{1}(t))|^2}\right|=
1+|g(\sigma_{1}(t))|\frac{|\alpha_{2}(t)|}{|1-|g(\sigma_{1}(t))|^2|}\leq (T+1).
\end{split}
\end{equation*}
and since $T\to \frac{T}{1+T}$ is a growing function, we  have
that
\[ \frac{|\alpha_{2}(t)|}{|1-\overline{g(\sigma_{1}(t))}\sigma_{2}(t)|}=
\frac{|\alpha_{2}(t)|}{1-|g(\sigma_{1}(t))|^2}\;
\frac{1-|g(\sigma_{1}(t))|^2}{|1-\overline{g(\sigma_{1}(t))}\sigma_{2}(t)|}\geq\frac{T}{1+T}>
\frac{\varepsilon}{1+\varepsilon}>0\] which contradicts the
hypothesis of $g$-speciality.
\end{proof}
We have now the following result of Lindel\"of type for Busemann
functions:
\begin{theorem}\label{LindeloffBusemannTheorem}
Let $f\in\Hol(\Delta^2,\Delta)$ be a holomorphic function. Given
$x\in\partial\Delta^2$ let $\varphi_{g}$ be a complex geodesic
passing through $x$ and $(\varphi_{g},\pi_{g})$ a $g-$projection device. Assume that  $f$ is
$K_{g}-$bounded. If $\sigma_{0}$ is a  $g$-special and $g$-restricted
 $x-$curve such that
$$\lim\limits_{t\to 1^-}f(\sigma_{0}(t))=L$$
then $f$ admits restricted $\tilde{K_{g}}-$limit equal to $L$ at $x.$
\end{theorem}
\begin{proof}
Let us consider  a $g$-special and
$g$-restricted
 $x-$curve $\sigma.$
 By definition there exists a constant $M>1$ such that $\tilde{\pi_{g}}(\sigma(t))$ approaches  $x_{1}$ inside
 a Stolz region $H(x_{1},M).$
We claim that
\begin{equation}\label{primopassoLindelof2}
\forall M_{1}>M,\:\:
K_{H_{\varphi_{g}}(x,M_{1})}(\sigma(t),\pi_{g}(\sigma(t)))\to 0
\:\:\hbox{as} \:\: t\to 1^-.
\end{equation}
For any $t\in[0,1)$ let us consider the map $\psi_{t}:\C \to\C^2$ given by
$$\phi_{t}(z)=\pi_{g}(\sigma(t))+z[\sigma(t)-\pi_{g}(\sigma(t))].$$
Let us notice that $\phi_{t}(0)=\pi_{g}(\sigma(t))$ and
$\phi_{t}(1)=\sigma(t).$ We claim  that the following statement is
true:
\begin{equation}\label{AffermazioneAdiLindelof2}
\begin{split}
\forall\;R>0\:\:\exists\:\: t_{0}=t_{0}(R)&\in [0,1) \:\:\hbox{such that} \:\: \forall \:t\in [0,1)\:\: :\:\: t>t_{0}(R)
\\
&\phi_{t}(\Delta_{R})\subset H_{\varphi_{g}}(x,M_{1}).
\end{split}
\end{equation}
Assuming (\ref{AffermazioneAdiLindelof2}) we get:
\begin{equation}\label{R(t)}
R(t):=\sup\{r>0:\:\varphi(\Delta_{r})\subset H_{\varphi_{g}}(x,M_{1})\}\to\infty
\end{equation}
as $t\to 1^-,$ and since, by the very definition \\
\begin{equation*}
\begin{split}
&K_{H_{\varphi_{g}}(x,M_{1})}(\sigma(t), \pi_{g}(\sigma(t))) \leq\\
&\leq \inf\{\frac{1}{R}\: :\: \exists
\varphi_{g}\in  \Hol(\Delta_{R},H_{\varphi_{g}}(x,M_{1}))\: :\: \varphi_{g}(0)=\pi_{g}(\sigma(t))\:\hbox{and}
\:\varphi_{g}(1)=\sigma(t)\}
\end{split}
\end{equation*}
then equation (\ref{primopassoLindelof2}) follows from equation
(\ref{R(t)}) and statement (\ref{AffermazioneAdiLindelof2}). Thus
we are left to prove (\ref{AffermazioneAdiLindelof2}). Assume  by
contradiction that  (\ref{AffermazioneAdiLindelof2}) is false.
Then there exist $M_{1}>M$ and $R_{0}>1$ such that for any
$t_{0}\in [0,1)$ there are $t'=t'(t_{0})\in(t_{0},1)$ and
$z_{0}=z_{0}(t_{0})\in\Delta_{R_{0}}$ such that
$\psi_{t'}(z_{0})\notin H_{\varphi_{g}}(x,M_{1}).$ Moreover, by
Proposition \ref{lemma1preparazioneLindeloff2}, $\pi_{g}(\sigma(t))\in
H_{\varphi_{g}}(x,M_{1})$ eventually, and in particular we can
choose $t'_{0}=t'_{0}(R_{0})\in(0,1)$ such that
$\pi_{g}(\sigma(t'))\in H_{\varphi_{g}}(x,M_{1})$ for all
$t_{0}>t'_{0}.$ Being $H_{\varphi_{g}}(x,M_{1})$ open we can also
assume that $\psi_{t'}(z_{0})\in\partial H_{\varphi_{g}}(x,M_{1})$
but $\psi_{t'}(z)\in H_{\varphi_{g}}(x,M_{1})$ for all $z\in
\Delta_{|z_{0}|}.$
\begin{remark}\label{successione non sul bordo}
Let us notice that there exists $t''_{0}>0$ such that
$\psi_{t'}(z)\notin\partial\Delta^2$  and
$\psi_{t'}(z)\in\Delta^2$ for all $t'_{0}>t''_{0}.$ Indeed,
suppose by contradiction that, for any $t_{0}\in [0,1)$ there are
$t'=t'(t_{0})\in(t_{0},1)$ and
$z_{0}=z_{0}(t_{0})\in\Delta_{R_{0}}$ such that
$\psi_{t'}(z_{0})\in \partial
H_{\varphi_{g}}(x,M_{1})\cap\partial\Delta^2.$ This implies that it
is possible to construct  two sequences, say
$\{t'_{k}\}_{k\in\N}=\{t'_{k}(t_{0})\}_{k\in\N}\subset(t_{0},1)$
and
$\{z^{k}_{0}\}_{k\in\N}=\{z^{k}_{0}(t_{0})\}_{k\in\N}\subset\Delta_{R_{0}}$
such that $\psi_{t'_{k}}(z^{k}_{0})\in \partial
H_{\varphi_{g}}(x,M_{1})\cap\partial\Delta^2.$ Since
$$\phi_{t'_k}(z^k_{0})=\pi_{g}(\sigma(t'_k))+z^k_{0}[\sigma(t'_k)-\pi_{g}(\sigma(t'_k))]=
(\sigma_{1}(t'_k),g(\sigma_{1}(t'_k))+z^k_{0}\alpha_{2}(t'_k))$$
it follows that $$|g(\sigma_{1}(t'_k))+\alpha_{2}(t'_k)|=1.$$ In
particular
\begin{equation*}
\begin{split}
0=&\frac{1-|g(\sigma_{1}(t'_k))+z^k_{0}\alpha_{2}(t'_k)|}{1-|g(\sigma_{1}(t'_k))|}\geq
\frac{1-|g(\sigma_{1}(t'_k))|-|z^k_{0}||\alpha_{2}(t'_k)|}{1-|g(\sigma_{1}(t'_k))|}\\=&
1-\frac{|z^k_{0}||\alpha_{2}(t'_k)|}{1-|g(\sigma_{1}(t'_k))|}
\geq 1-\frac{R_{0}|\alpha_{2}(t'_k)|}{1-|g(\sigma_{1}(t'_k))|}\to 1^-
\end{split}
\end{equation*}
as $t\to 1^-,$  a contradiction.
\end{remark}
According to remark \ref{successione non sul bordo} and by definition of g-Koranyi region,
 we can write
 \begin{equation}\label{M1}
\log M_{1}=\lim\limits_{s\to 1^-}K_{\Delta^2}(\psi_{t'}(z_{0}),\varphi(s))-\omega(0,s)+
K_{\Delta^2}(\psi_{t'}(z_{0}),\varphi(0)).\end{equation}
Furthermore for any $z\in\Delta_{R_{0}}$
\begin{equation}\label{maggiorazioni per la distanza di log M1}
\begin{split}
&K_{\Delta^2}(\psi_{t'}(z),\varphi(s))-\omega(0,s)+
K_{\Delta^2}(\psi_{t'}(z),\varphi(0))=
\\ &=
\max\{\omega(\sigma_{1}(t'),s); \omega(g(\sigma_{1}(t'))+\alpha_{2}(t')z,g(s))\}
-\omega(0,s)+\\&\hspace{4.5cm}+ \max\{\omega(\sigma_{1}(t'),0); \omega(g(\sigma_{1}(t'))+\alpha_{2}(t')z,g(0))\}\leq\\
&\leq \omega(g(\sigma_{1}(t'))+\alpha_{2}(t')z,g(\sigma_{1}(t'))+\omega(\sigma_{1}(t'),s)-\omega(0,s)+
\\&\hspace{4.5cm}+ \omega(\sigma_{1}(t'),0)+\omega(g(\sigma_{1}(t'))+\alpha_{2}(t')z,g(\sigma_{1}(t')))\}=\\
&=\omega(\sigma_{1}(t'),s)-\omega(0,s)+ \omega(\sigma_{1}(t'),0)+2\omega(g(\sigma_{1}(t'))
+\alpha_{2}(t')z,g(\sigma_{1}(t'))).
\end{split}
\end{equation}
Let us observe that
\begin{equation}\label{distanza a 0 in Lindeloff 2}
\lim\limits_{t'\to 1^-}\omega(g(\sigma_{1}(t'))
+\alpha_{2}(t')z,g(\sigma_{1}(t')))=0
\end{equation}
uniformly for $z\in\Delta_{R_{0}}.$
Indeed by definition of Poincar\'{e} distance, we have
\begin{equation*}
\begin{split}
\lim\limits_{t'\to 1^-}\omega(g(\sigma_{1}(t'))
+\alpha_{2}(t')z,g(\sigma_{1}(t')))=\lim\limits_{t'\to 1^-}\frac{1}{2}\log\frac{1+\left|\frac{\alpha_{2}(t')z}
{1-\overline{g(\sigma_{1}(t'))}(g(\sigma_{1}(t'))+\alpha_{2}(t')z)}\right|}{1-\left|\frac{\alpha_{2}(t')z}
{1-\overline{g(\sigma_{1}(t'))}(g(\sigma_{1}(t'))+\alpha_{2}(t')z)}\right|}
\end{split}
\end{equation*}
and the argument of this logarithm tends to $1$ since
\begin{equation*}
\begin{split}
&\left|\frac{\alpha_{2}(t')z}
{1-\overline{g(\sigma_{1}(t'))}(g(\sigma_{1}(t'))+\alpha_{2}(t')z)}\right|
\leq \frac{|\alpha_{2}(t')z|}{1-|g(\sigma_{1}(t'))||(g(\sigma_{1}(t'))+\alpha_{2}(t')z)|}\\
& \leq \frac{|\alpha_{2}(t')z|}{1-|(g(\sigma_{1}(t'))+\alpha_{2}(t')z)|}
\leq \frac{|\alpha_{2}(t')z|}{1-|g(\sigma_{1}(t'))|-|\alpha_{2}(t')z|}\\
&\leq \frac{1}{\frac{1-|g(\sigma_{1}(t'))|}{|\alpha_{2}(t')z|}-1}\to 0 \:\hbox{as}\:t\to 1^-
\end{split}
\end{equation*}
by Lemma \ref{lemma2preparazioneLindeoff2}.
Thus equation (\ref{distanza a 0 in Lindeloff 2}) is proved. In particular it is true for $z=z_{0}$
and then by equations (\ref{M1}) and  (\ref{maggiorazioni per la distanza di log M1}) we get
\[
\log M_{1}\leq \lim\limits_{s\to 1^-} \omega(\sigma_{1}(t'),s)-\omega(0,s)+ \omega(\sigma_{1}(t'),0)+2\omega(g(\sigma_{1}(t'))
+\alpha_{2}(t')z,g(\sigma_{1}(t')))  \]
and in particular, eventually,
\[M<M_{1}\leq \frac{|1-\sigma_{1}(t')|}{1-|\sigma_{1}(t')|}\]
but it is a contradiction since $\sigma$
is $g$-restricted.
This concludes the proof of  (\ref{AffermazioneAdiLindelof2}).Now, since $f$ is a
$K_{g}-$bounded function, there exists $c>0$ such that
\[K_{\Delta_{c(M_{1})}}(f(\sigma_{0}(t)),f(\pi_{g}(\sigma_{0}(t))))
\leq K_{H_{\varphi_{g}}}(\sigma_{0}(t),\pi_{g}(\sigma_{0}(t)))\to 0 \:\hbox{as}\: t\to 1^-.\]
Now we can proceed as in Theorem \ref{Lindelof1}
to complete the proof.
\end{proof}


\section{Julia's Lemma}
We want to give a new generalization to polydiscs of the Julia's lemma, using the Busemann functions.
The idea is to consider the rate of approach
of $f$ along particular directions given by geodesics passing through $x.$
\begin{definition}\label{definizione lambdaj}
Let $f\in\Hol(\Delta^2,\Delta)$ and  $x\in\partial\Delta^2.$ Let us consider
a complex geodesic $\varphi_{g}\in\Hol(\Delta,\Delta^2)$   passing through $x.$
Let $\lambda_{g}$ be the boundary dilation coefficient of $g$ at $x_{1}.$
The number $\lambda_{\varphi_{g}}(f)$ defined by
$$\frac{1}{2}\log\lambda_{\varphi_{g}}(f):=
\lim\limits_{t\to 1^-} K_{\Delta^2}(0,\varphi_{g}(tx_{1}))-\omega(0,f(\varphi_{g}(tx_{1}))).$$
 is the \emph{$\varphi_{g}-$boundary dilation coefficient} of $f$ at $x.$
\end{definition}
First we show  that $\lambda_{\varphi_{g}}(f)$ is well defined.
We prove this fact studying separately two cases:

$1)$ $x\in(\partial\Delta)^2$ or

$2)$ $x\in[(\partial\Delta^2)-(\partial\Delta)^2].$

In the first case we can assume $x=(1,1)$ and in the second one we suppose $x=(1,0).$
\begin{remark}\label{coeff. in una variabile lungo t}
Suppose (as in case $1)$) that $g\in\Hol(\Delta,\Delta)$ has non-tangential limit $1$ at
the point $1.$ Let observe that if $\lambda_{g}<\infty$ then
$\lambda_{g}=\lim_{t\to 1}\frac{1-|g(t)|}{1-t}.$ Indeed  we have
\[ \lambda_{g}:=\liminf_{z\to 1}\frac{1-|g(z)|}{1-|z|}\leq
 \liminf_{t\to 1}\frac{1-|g(t)|}{1-t}\leq \limsup_{t\to 1}\frac{1-|g(t)|}{1-t}\]
 and by the triangular inequality we get
 \[\leq \limsup_{t\to 1}\frac{1-|g(t)|}{1-t}
 \leq \lim_{t\to 1}\frac{|1-g(t)|}{1-t}=\lambda_{g}\]
by the classical Julia-Wolff-Carath\'{e}odory theorem.
\end{remark}

Let us consider the case $1)$ and let $\psi_{g}(z)=(\theta(z),g(\theta(z)))$ be
 another parametrization of  the geodesic $\varphi_{g},$ with $\theta\in\Aut(\Delta).$
We notice that $\theta(1)=1.$
As a matter of notation we call respectively $\lambda_{\varphi}$ and $\lambda_{\psi}$ the boundary dilation coefficient
of $f$ at $x$ computed with respect to the parameterizations  $z\to (z,g(z))$ and
$z\to (\theta(z),g(\theta(z))).$
By the above definition we have that:
\begin{equation}\label{limiteBuonaDefCoeffDilatazione}
\begin{split}
\frac{1}{2}\log\lambda_{\varphi}&=\lim\limits_{t\to 1^-}[K_{\Delta^2}(0,\varphi_{g}(t))-\omega(0,f(\varphi_{g}(t)))]
                              \\&=\lim\limits_{t\to 1^-}[K_{\Delta^2}(0,\psi_{g}(\theta^{-1}(t)))
                              -\omega(0,f(\psi_{g}(\theta^{-1}(t)))]
                              \\&=\lim\limits_{t\to 1^-}[\max\{\omega(0,t), \omega(0,g(t))\}
                              -\omega(0,f(\psi_{g}(\theta^{-1}(t)))].
\end{split}
\end{equation}
If $\lambda_{g}\geq 1,$ then $\max\{\omega(0,t), \omega(0,g(t))\}=\omega(0,t)$ and the last member of equation
(\ref{limiteBuonaDefCoeffDilatazione}) becomes
\begin{equation}\label{calcolo esplicito coeff. con parametrizzazione}
\begin{split}
&\lim\limits_{t\to 1^-}[\omega(0,t)-\omega(0,\theta^{-1}(t))+ \omega(0,\theta^{-1}(t))
-\omega(0,f(\psi_{g}(\theta^{-1}(t))))]
\\&
\geq  \liminf\limits_{t\to 1^-}[\omega(0,t)-\omega(0,\theta^{-1}(t))]+
\liminf\limits_{z\to 1^-}[\omega(0,z)
-\omega(0,f(\psi_{g}(z)))]
\\&
\geq  \lim\limits_{t\to 1^-}[\omega(0,t)-\omega(0,\theta^{-1}(t))]+
\lim\limits_{t\to 1^-}[\omega(0,t)-\omega(0,f(\psi_{g}(t)))]
\end{split}
\end{equation}
Since $f\circ\psi_{g}$ is a holomorphic self map of the unit disc,
by remark \ref{coeff. in una variabile lungo t}, the equation
(\ref{calcolo esplicito coeff. con parametrizzazione}) becomes
\begin{equation}
\begin{split}
&\lim\limits_{t\to 1^-}[\omega(0,t)-\omega(0,\theta^{-1}(t))]+
\lim\limits_{t\to 1^-}[\omega(0,t)-\omega(0,f(\psi_{g}(t)))]
\\
& =\lim\limits_{t\to 1^-}[\omega(0,t)-\omega(0,\theta^{-1}(t))]+
\lim\limits_{t\to 1^-}[\omega(0,t)-K_{\Delta^2}(0,\psi_{g}(t))]
\\&\hspace{5cm}
+\lim\limits_{t\to
1^-}[K_{\Delta^2}(0,\psi_{g}(t))-\omega(0,f(\psi_{g}(t)))]
\\&=\lim\limits_{t\to 1^-}[\omega(0,t)-\omega(0,\theta^{-1}(t))]+
\lim\limits_{t\to 1^-}[\omega(0,t)-\omega(0,\theta(t))]
\\&\hspace{5cm}
+\lim\limits_{t\to
1^-}[K_{\Delta^2}(0,\psi_{g}(t))-\omega(0,f(\psi_{g}(t)))]
\\& =\frac{1}{2}\log(\lambda_{\theta^{-1}}\:\lambda_{\theta}\:\lambda_{\psi})
=\frac{1}{2}\log(\lambda_{\psi})
\end{split}
\end{equation}

If $\lambda_{g}\leq 1,$ then $\max\{\omega(0,t), \omega(0,g(t))\}=\omega(0,g(t))$ and the last member of equation
(\ref{limiteBuonaDefCoeffDilatazione}) becomes
\begin{equation*}
\begin{split}
&\lim\limits_{t\to 1^-}\omega(0,g(t))-\omega(0,t)+\omega(0,t)-\omega(0,\theta^{-1}(t))+ \omega(0,\theta^{-1}(t))
-\omega(0,f(\psi_{g}(\theta^{-1}(t)))
\\&
=\frac{1}{2}\log\frac{1}{\lambda_{g}\lambda_{\theta}}+\liminf\limits_{z\to 1^-}\omega(0,z)
-\omega(0,f(\psi_{g}(z))
\\&
=\frac{1}{2}\log\frac{1}{\lambda_{g}\lambda_{\theta}}+\lim\limits_{t\to
1^-}\omega(0,t) -K_{\Delta^2}(0,\psi_{g}(t))
+K_{\Delta^2}(0,\psi_{g}(t))-\omega(0,f(\psi_{g}(t))
\\&
=\frac{1}{2}\log(\frac{1}{\lambda_{\theta}\lambda_{g}}\:\lambda_{\theta}\lambda_{g}\:\lambda_{\psi})
=\frac{1}{2}\log(\lambda_{\psi})
\end{split}
\end{equation*}
Then we have $\lambda_{\varphi}\geq \lambda_{\psi}.$
Swapping the roles of $\lambda_{\varphi}$ and $\lambda_{\psi}$ in the above inequalities, we also get that
$\lambda_{\psi}\geq \lambda_{\varphi}$ and thus $\lambda_{\varphi}= \lambda_{\psi}.$

In case $2)$ we can suppose $x=(1,0),$  and we consider a complex geodesic
$\varphi_{g}$ passing through $(1,0)$ parameterized by $\varphi_{g}(z)=(z,g(z)).$
  We also consider
another parametrization $\psi_{g}(z)=(\theta(z),g(\theta(z))),$
with $\theta(z)\in \Aut(\Delta).$
We notice that $\theta(1)=1.$
Since, in this case, $\lambda_{g}=\infty,$
we can repeat the calculation done in (\ref{limiteBuonaDefCoeffDilatazione}), and in
(\ref{calcolo esplicito coeff. con parametrizzazione}) obtaining
$\lambda_{\psi}= \lambda_{\varphi}.$
Thus $\lambda_{\varphi_{g}}(f)$ is well defined.
Furthermore we have an interesting property.
Let be $\frac{1}{2}\log \alpha_{f}:=\liminf\limits_{w\to x}[K_{\Delta^n}(0,w)-\omega(0,f(w))].$
Let notice that $\alpha(f)$ is  the boundary dilation coefficient of $f$ at $x,$ defined by Abate in
\cite{ArticoloJWCAbate} and the following property holds (see \cite{ArticoloJWCAbate} for the proof)
$$\frac{1}{2}\log\alpha(f)=\lim\limits_{t\to 1}[K_{\Delta^n}(0,tx)-\omega(0,f(tx))].$$

\begin{theorem}\label{JuliaBusemannlemma}
Let $f=(f_{1},f_{2})\in\Hol(\Delta^2,\Delta^2)$. Let $x=(x_{1},
x_{2})\in\partial(\Delta\times\Delta)=\partial\Delta^2$
and let (for example) $\varphi_{g}=(z,g(z))$ be a complex
geodesic passing through $x.$
Let
$$\frac{1}{2}\log\lambda_{j}:=\lim\limits_{t\to
x}[K_{\Delta^2}(0,\varphi_{g}(tx_{1}))
-\omega(0,f_{j}(\varphi_{g}(tx_{1}))]\:\:\: j=1,2$$
Suppose that either $\lambda_{1}<\infty$ or $\lambda_{2}<\infty.$ Then
there exists  a point $y=(y_{1},y_{2})\in \partial\Delta^2$  such
that  for all $R>0$
$$f(\B_{(1,\lambda_{g})}(x,R))\subseteq\B_{(\lambda_{1},\lambda_{2})}(y,R).$$
\end{theorem}
%

\begin{proof}
Let us first suppose  that $\lambda_{j}<\infty,$ for $j=1,2$ then
$$\frac{1}{2}\log\alpha_{j}=\liminf\limits_{z\to x}K_{\Delta^2}(0,z)-\omega(0,f_{j}(z))<\infty\:\:\: j=1,2.$$
As shown by Abate in theorem $3.1$ in (\cite{ArticoloJWCAbate}), we can choose a sequence $z_{\nu}\in \Delta^{2},$
converging to $x,$ such that
$$\lim\limits_{\nu\to \infty}K_{\Delta^2}(0,z_{\nu})-\omega(0,f_{j}(z_{\nu}))=
\liminf\limits_{z\to x}K_{\Delta^2}(0,z)-\omega(0,f_{j}(z)).$$
Up to a subsequence, we can assume that $f_{j}(z_{\nu})\to y_{j}\in\overline\Delta.$
Since $\Delta^2$ is complete hyperbolic, we have that $K_{\Delta^2}(0,z_{\nu})\to +\infty;$
therefore $\omega(0,f_{j}(z_{\nu}))\to +\infty$ as well, and $y_{j}\in\partial\Delta.$
Thus there exists a point $y=(y_{1},y_{2})\in(\partial\Delta)^2$ such that $f_{j}(z)\to y_{j}$ as $z\to 1,$ $j=1,2.$

We claim that
$$f(\B_{(1,\lambda_{g})}(x,R))\subseteq\B_{(\lambda_{1},\lambda_{2})}(y,R)\:\:\forall\:R>0.$$
Without loss of generality let us suppose that $x_{1}=1.$ Fix $z\in\B_{(1,\lambda_{g})}(x,R)).$
We have, for $j=1,2,$ that
\begin{equation*}
\begin{split}
&\lim\limits_{w\to y_{j}}\omega(f_{j}(z),w)-\omega(0,w)=\lim\limits_{s\to 1}\omega(f_{j}(z),f_{j}
(\varphi_{g}(s))-\omega(0,f_{j}(\varphi_{g}(s))) \\
&\leq \liminf\limits_{s\to 1}K_{\Delta^2}(z,\varphi_{g}(s))-\omega(0,f_{j}(\varphi_{g}(s)))\\&
=\liminf\limits_{s\to 1}K_{\Delta^2}(z,\varphi_{g}(s))-\omega(0,s)+\omega(0,s)
-K_{\Delta^2}(0,\varphi_{g}(s))+\\&\hspace{7cm}+K_{\Delta^2}(0,\varphi_{g}(s))
-\omega(0,f_{j}(\varphi_{g}(s)))\\&
\leq\lim\limits_{s\to 1}K_{\Delta^2}(z,\varphi_{g}(s))-\omega(0,s)+\lim\limits_{s\to 1}K_{\Delta^2}(0,\varphi_{g}(s))
-\omega(0,f_{j}(\varphi_{g}(s)))\leq \frac{1}{2}\log\lambda_{j}R.
\end{split}
\end{equation*}
Then $ \forall\:R>0$
$$f(\B_{(1,\lambda_{g})}(x,R))=f(E(x_{1},R)\times E(x_{2},\lambda_{g}R))\subseteq E(y_{1},\lambda_{1}R)\times
E(y_{2},\lambda_{2}R).$$

Let suppose now that $\lambda_{1}=\infty$ and $\lambda_{2}<\infty.$
By the above calculation we get that
$ \forall\:R>0$
$$f(\B_{(1,\lambda_{g})}(x,R))=f(E(x_{1},R)\times E(x_{2},\lambda_{g}R))\subseteq \Delta\times
E(y_{2},\lambda_{2}R).$$

\end{proof}
Let us notice that the following proposition holds:
\begin{proposition}\label{equivalenza coeff. dilatazione}
For all complex geodesic $\varphi_{g}$ passing through $x$ such that the coefficient
$\lambda_{g}<\infty$ we have that  $\lambda_{\varphi_{g}}(f)$ is finite if and only if $\alpha(f)$ is finite.
\end{proposition}
\begin{proof}
It is clear that $\lambda_{\varphi_{g}}(f)\geq \alpha(f)$ and thus if $\lambda_{\varphi_{g}}(f)$ is finite then
also $\alpha(f)$ does.
On the other hand let us suppose that $\alpha(f)$ is finite and let us
denote by $\varphi_{x}(z)$ the complex geodesic passing through the point $x=(x_{1},x_{2})\in\partial\Delta^2$
and parameterized by $z\to zx.$
Let us consider $\varphi_{g}(z)=(z,g(z))$ another complex geodesic passing through the point $x,$
and $\pi_{g}\in\Hol(\Delta^2,\Delta^2)$  the projection on the complex geodesic $\varphi_{g}(z)$ given by
$\pi_{g}(z_{1},z_{2})=(z_{1},g(z_{1})).$ Note that $\pi_{g}(\varphi_{x}(t))=\pi_{g}(tx)=(tx_{1},g(tx_{1}))
=\varphi_{g}(tx_{1}).$ Let us suppose, without loss of generality  that $x_{1}=1.$
Then
\begin{equation}
\begin{split}
\frac{1}{2}\log\lambda_{\varphi_{g}}(f)&=
\lim\limits_{t \to 1^-}K_{\Delta^2}(0,\varphi_{g}(t))-\omega(0,f(\varphi_{g}(t)))\\&=
\lim\limits_{t \to 1^-}K_{\Delta^2}(0,\pi_{g}(\varphi_{x}(t)))-\omega(0,f(\pi_{g}(\varphi_{x}(t))))\\
&=\lim\limits_{t \to 1^-}K_{\Delta^2}(0,\pi_{g}(\varphi_{x}(t)))-K_{\Delta^2}(0,\varphi_{x}(t))+
K_{\Delta^2}(0,\varphi_{x}(t))\\
&\hspace{1cm}-\omega(0,f(\varphi_{x}(t)))+\omega(0,f(\varphi_{x}(t)))-\omega(0,f(\pi_{g}(\varphi_{x}(t))))\\
&\leq \lim\limits_{t \to 1^-}K_{\Delta^2}(0,\varphi_{x}(t))-\omega(0,f(\varphi_{x}(t)))
\\&\hspace{5cm}+2K_{\Delta^2}(\varphi_{x}(t),\pi_{g}(\varphi_{x}(t)))\\
&=\lim\limits_{t \to 1^-}K_{\Delta^2}(0,\varphi_{x}(t))
-\omega(0,f(\varphi_{x}(t)))+2\omega(t,g(t))\\&=\frac{1}{2}\log\alpha(f)+\log
\lambda_{g}
\end{split}
\end{equation}
and in our setting $\lambda_{g}$ is finite and then if $\alpha(f)$ is finite also $\lambda_{\varphi_{g}}(f)$ does.
\end{proof}

\section{The Julia-Wolff-Carath\'{e}odory theorem}

We are finally ready to state and prove our generalization of the
 Julia-Wolff-Carath\'{e}odory theorem obtained using Busemann functions.
\begin{theorem}\label{JWCBusemannTheorem}
Let $f\in\Hol(\Delta^2,\Delta^2)$  and $x\in\partial\Delta^2.$
Let  $\varphi_{g}$ be any complex geodesic  passing through
$x$ and parameterized by $\varphi_{g}(z)=(z,g(z)),$ with $g\in\Hol(\Delta,\Delta)$
such that
$$\frac{1}{2}\log\lambda_{j}=\lim\limits_{t\to 1}K_{\Delta^2}(0,\varphi_{g}(tx_{1}))-\omega(0,f_{j}(\varphi_{g}(tx_{1})))<\infty.$$
for $j=1,2.$ Let
 $\tilde{\pi_{g}}(z):\Delta^2\to\Delta$ be the
 $g-$left-inverse
of $\varphi_{g}$ given by $\tilde{\pi_{g}}(z_{1},z_{2})=z_{1}.$
Then there exists a point $y=(y_{1}, y_{2})\in\partial\Delta^2$ such that
$$\displaystyle\tilde{K_{g}}-\!\lim_{\hspace{-.5cm}z\to x}\frac{y_{j}-f_{j}(z)}{1-\tilde{\pi}(z)}
=\lambda_{j}\min\{1,\lambda_{g}\}$$
$$\displaystyle\tilde{K_{g}}-\!\lim_{\hspace{-.5cm}z\to x}\frac{y_{j}-f_{j}(z)}{1-z_{2}}=\frac{\lambda_{j}}{\max\{1,\lambda_{g}\}}$$
\end{theorem}

To prove this theorem we  need first the following two lemmas:
\begin{lemma}\label{penultimo lemma}
Let $f\in\Hol(\Delta^2,\Delta^2)$  and $x\in\partial\Delta^2.$ Suppose $|x_{1}|=1.$
Suppose  there exists a complex geodesic $\varphi_{g}$ passing through
$x$ and parameterized by $\varphi_{g}(z)=(z,g(z)),$ with $g\in\Hol(\Delta,\Delta)$ such that
$$\frac{1}{2}\log\lambda_{j}=\lim\limits_{t\to 1}K_{\Delta^2}(0,\varphi_{g}(tx_{1}))-\omega(0,f_{j}(\varphi_{g}(tx_{1})))<\infty.$$
Let $\tilde{\pi_{g}}(z):\Delta^2\to\Delta$ be the
 $g-$left-inverse
of $\varphi_{g}$ given by $\tilde{\pi_{g}}(z_{1},z_{2})=z_{1}.$
Then there exists a point $y=(y_{1}, y_{2})\in\partial\Delta^2$ and   a constant,
say $c_{g}>0,$ depending on $g,$ such that, given $M>1,$ for all $z\in H_{\varphi_{g}}(x,M)$
$$\left|\frac{y_{j}-f_{j}(z)}{1-\tilde{\pi}(z)} \right|\leq 2\lambda_{1}M^2c_{g}\hspace{1cm}\hbox{and}
\hspace{1cm}\left|\frac{y_{j}-f_{j}(z)}{1-z_{2}}\right|\leq 2\lambda_{1}M^2c_{g} .$$
\end{lemma}
\begin{proof}
Without loss of generality let us suppose that $x_{1}=1.$
Let $z\in H_{\varphi_{g}}(x,M)$ and set $\frac{1}{2}\log R:=\log M-K_{\Delta^2}(\varphi_{g}(0),z).$
Thus
$$\lim\limits_{s\to 1}K_{\Delta^2}(z,\varphi_{g}(s))-K_{\Delta^2}(\varphi_{g}(0),\varphi_{g}(s))<\frac{1}{2}\log R$$
which implies $z\in\B_{(1,\lambda_{g})}(x,R).$
By Lemma \ref{JuliaBusemannlemma}  there exists a point $y=(y_{1}, y_{2})\in\partial\Delta^2$ and a complex geodesic $\varphi_{\tilde{g}}$ passing through
$y$ and parameterized by $\varphi_{\tilde{g}}(z)=(z,\tilde{g}(z)),$ with $\tilde{g}\in\Hol(\Delta,\Delta)$ such that
$\lambda_{\tilde{g}}=\frac{\lambda_{2}}{\lambda_{1}}$ and
$f(z)\in\B_{(1,\lambda_{\tilde{g}})}(y,\lambda_{1}R).$ Without loss of generality let us suppose that $y=(1,1).$
In particular, by the very definition of Busemann sublevel sets, we have
\[ \frac{1}{2}\log \lambda_{1}R\geq \lim\limits_{s\to 1}K_{\Delta^2}(f(z),\varphi_{\tilde{g}}(s))-\omega(0,s)\]
\[\geq \lim\limits_{s\to 1}\omega(f_{j}(z),\varphi_{\tilde{g}_{j}}(s))-\omega(0,s) \]
for $j=1,2.$ Moreover let us notice that
\[-\omega(0,f_{j}(z))\leq \omega(f_{j}(z),s)-\omega(0,s)\:\:\forall\;s\in(0,1).\]
For sake of clearness we argue for $j=1$
\[\lim\limits_{s\to 1}\omega(f_{1}(z),s)-\omega(0,s)-\omega(0,f_{1}(z))\leq \log \lambda_{1}R\]
then
\[\frac{|1-f_{1}(z)|^2}{1-|f_{1}(z)|^2}\;\frac{1-|f_{1}(z)|}{1+|f_{1}(z)|}=\left[\frac{|1-f_{1}(z)|}{1+|f_{1}(z)|}
\right]^2\leq(\lambda_{1}R)^2.\]
 Furthermore we know that
 \[-2 K_{\Delta^2}(\varphi_{g}(0),z)\leq 2K_{\Delta^2}(\varphi_{g}(0),0)-2K_{\Delta^2}(0,z)\]
 \[=
 \log\frac{1+|||\varphi_{g}(0)|||}{1-|||\varphi_{g}(0)|||}\;\frac{1-|||z|||}{1+|||z|||}=\log\frac{1+|g(0)|}
 {1-|g(0)|}\;\frac{1-|||z|||}{1+|||z|||}.\]
 and , by definition of $R,$ set $c_{g}:=\frac{1+|g(0)|}
 {1-|g(0)|}$ then we get
\[ \frac{|1-f_{1}(z)|}{1+|f_{1}(z)|}\leq \lambda_{1}M^2c_{g}\frac{1-|||z|||}{1+|||z|||}\leq
\lambda_{1}M^2c_{g}\frac{1-|z_{i}|}{1+|z_{i}|}  \:\:\hbox{for}\:\: i=1,2
\]
If $i=1,$ being  $\tilde{\pi}(z_{1},z_{2})=z_{1},$ then
\[ \frac{|1-f_{1}(z)|}{|1-\tilde{\pi}(z)|}\leq \frac{|1-f_{1}(z)|}{1-|\tilde{\pi}(z)|}\leq
\lambda_{1}M^2c_{g}\frac{1+|f_{1}(z)|}{1+|z_{i}|}\leq 2\lambda_{1}M^2c_{g}\]
and if $i=2$
\[ \frac{|1-f_{1}(z)|}{|1-z_{2}|}\leq \frac{|1-f_{1}(z)|}{1-|z_{2}|}\leq 2\lambda_{1}M^2c_{g}.\]
With the same techniques we proved the statement for the second component $f_{2}.$

\end{proof}

\begin{lemma}\label{ultimo}
Let $f\in\Hol(\Delta^2,\Delta^2)$ be a holomorphic function and $x\in\partial\Delta^2.$
Suppose there exists a complex geodesic $\varphi_{g}$ passing through
$x$ and parameterized by $\varphi_{g}(z)=(z,g(z)),$ with $g\in\Hol(\Delta,\Delta)$ such that
$$\frac{1}{2}\log\lambda_{j}=\lim\limits_{t\to 1}K_{\Delta^2}(0,\varphi_{g}(tx_{1}))-\omega(0,f_{j}(\varphi_{g}(tx_{1})))<\infty.$$
Let $\tilde{\pi_{g}}(z):\Delta^2\to\Delta$ be a $g-$left inverse
of $\varphi_{g}$ given by $\tilde{\pi_{g}}(z_{1},z_{2})=z_{1}.$
Then there exists a point $y=(y_{1}, y_{2})\in\partial\Delta^2$ such that, for $j=1,2,$
\begin{equation}\label{limite1ultimolemma}\lim\limits_{s\to 1}\frac{1-f_{j}(\varphi_{g}(sx_{1}))}{1-sx_{1}}=\lambda_{j}\min\{1,\lambda_{g}\}\end{equation}
\begin{equation}\label{limite2ultimolemma}\lim\limits_{s\to 1}\frac{1-f_{j}(\varphi_{g}(sx_{1}))}{1-g(sx_{1})}=\frac{\lambda_{j}}{\max\{1,\lambda_{g}\}}.\end{equation}

\end{lemma}
\begin{proof}
Let us suppose that $x_{1}=1.$
Let $y$ the point given in Theorem \ref{JuliaBusemannlemma} and without loss of generality let us suppose that $y=(1,1)$.
By definition of $\lambda_{j}$ we have that
$$\lambda_{j}=\lim\limits_{s\to1}\frac{1-|f_{j}(\varphi_{g}(s))|}{1-|||\varphi_{g}(s)|||}.$$
Moreover the  limit (\ref{limite1ultimolemma}) exists for the
classic Julia-Wolff-Carath\'{e}odory theorem  and also the limit
(\ref{limite2ultimolemma}) exists, since
$$\frac{1-f_{j}(\varphi_{g}(s))}{1-g(s)}={\frac{1-f_{j}(\varphi_{g}(s))}{1-s}}\;\frac{1-s}{1-g(s)}.$$
Let proceed considering the following different cases $(a)$
 $\lambda_{g}\geq 1$ and
$(b)$ $\lambda_{g}\leq 1.$ Let study the case $(a).$
In this setting there exists a sequence $\{s_k\}_{k\in\N}\in(0,1)$ such that $s_k\to1$ as $k\to \infty$
and $\lim_{k\to\infty}[\omega(0,s_k)-\omega(0,g(s_k))]\geq 0$
which implies that $|||\varphi_{g}(s_k)|||=s_k$ and then
\[\lambda_{j}\min\{1,\lambda_{g}\}=\lim\limits_{s\to 1}\frac{1-|f_{j}(\varphi_{g}(s))|}{1-|||\varphi_{g}(s)|||}=
\lim\limits_{k\to\infty}\frac{1-|f_{j}(\varphi_{g}(s_{k}))|}{1-|||\varphi_{g}(s_{k})|||}\]\[=
\lim\limits_{k\to\infty}\frac{1-|f_{j}(\varphi_{g}(s_{k}))|}{1-s_{k}}=
\lim\limits_{s\to 1}\frac{1-|f_{j}(\varphi_{g}(s))|}{1-s}\]
and moreover,
\[\frac{1-f_{j}(\varphi_{g}(s))}{1-g(s)}={\frac{1-f_{j}(\varphi_{g}(s))}{1-s}}\;\frac{1-s}{1-g(s)}=\frac{\lambda_{j}}
{\lambda{g}}=\frac{\lambda_{j}}{\max\{1,\lambda_{g}\}}. \]
Using the same techniques we proved that the above equalities hold also in case $(b).$
\end{proof}
And now we are ready to prove Theorem \ref{JWCBusemannTheorem}.
\begin{proof}[Proof of Theorem \ref{JWCBusemannTheorem}]
Let us suppose $x_{1}=1.$ By Lemma \ref{ultimo} we have
\[\lim\limits_{s\to 1}\frac{1-|f_{j}(\varphi_{g}(s))|}{1-s}=\lambda_{j}\min\{1,\lambda_{g}\}\]
and   by Lemma \ref{penultimo lemma} we know that the function $\left|\frac{1-f_{j}(z)}
{1-\tilde{\pi}(z)}\right|$ is $K_{g}-$bounded.
Then, since the curve $\varphi_{g}(s)$ is $g$-special and $g$-restricted,
the conclusion of the proof follows by theorem \ref{LindeloffBusemannTheorem} .
\end{proof}


\section{Application to the dynamics}

Let $f\in\Hol(\Delta,\Delta)$ be without fixed points in $\Delta.$
The classical Wolff lemma ensures the existence of a unique point
$\tau\in\partial\Delta$ such that every horocycle centered in
$\tau$ is sent in itself by $f.$ The point $\tau$ is called the
\emph{Wolff point} of $f.$

Let $n\in\N,$ and set $f^n=f\circ\cdot\cdot\cdot\circ f$ the composition of $f$ with itself $n-$times.
We say that $\{f^n\}_{n\in\N}$ is the \emph{sequence of iterates} of $f.$
The Wolff-Denjoy lemma says that $\{f^n\}$ converges uniformly on compacta to the Wolff point $\tau.$

If we call \emph{target set}, $T(f),$ the set of the limit points of the sequence of the iterates and we denote
by $W(f)$ the set of the Wolff points of $f,$ then in one complex variable, we have that $T(f)\equiv W(f)=\{\tau\}.$

In \cite{Frosini2} we considered $f\in\Hol(\Delta^2,\Delta^2)$ without fixed points in $\Delta^2$ and we
defined the Wolff points of $f$ using the small and big horospheres.
\begin{definition}\cite{Frosini2}
Let $f\in\Hol(\Delta^2,\Delta^2)$ be without fixed points in $\Delta^2.$ A point $\tau\in\partial\Delta^2$
is a \emph{Wolff point} for $f$ if
$f(E(\tau,R))\subseteq E(\tau,R),$ for all $ R>0.$
\end{definition}
In this setting, in \cite{Frosini2} (see also \cite{Frosini}) we characterized the set of the Wolff points,
$W(f),$ for a holomorphic self map $f$ of the bidisc without fixed points.
As a spinning result (see also Herv\'{e} \cite{HerveSuBidisco}), we get that $W(f)\subset T(f),$
where $T(f)$ is the
\emph{target set} of $f$ defined as follows (\cite{Frosini},\cite{Frosini2}):
\begin{center}$T(f):= \{ x \in\overline{\Delta^{2}}  : \exists \; \{ k_{n}\}\subset \N
,\;\exists\; z \in \; \Delta^{2} \; \hbox{such that} \; f^{k_{n}}(z)\to x \;
\hbox{as}
  \; n\to \infty\}.$
  \end{center}
It turns out that this result can be improved using the Busemann functions.
\begin{definition}
Let $\tau\in \partial \Delta^2.$ We say that $\tau$ is a \emph{generalized Wolff point} for $f$
if there exists a geodesic $\varphi_{g},$ passing through the point $\tau,$ such that
every Busemann sublevel set, $\B_{(1,\lambda_{g})}(\tau,R),$
is sent in itself by $f,$ that is
$f(\B_{(1,\lambda_{g})}(\tau,R))\subseteq \B_{(1,\lambda_{g})}(\tau,R)$ for every $R>0.$

Let denote by $W_{G}(f)$ the set of the \emph{generalized Wolff points} for $f.$
\end{definition}

\begin{remark}\label{Nota su Wolff e Generalized}
Let notice that if $\tau$ is contained in a flat component of the boundary then
 $\tau$ is a Wolff point if and only if $\tau$ is a generalized Wolff point.
On the other hand, let consider a point  $\tau$  of the S\v{i}lov boundary of the bidisc.
If $\tau$ is a Wolff point for $f$ then i is also a generalized Wolff point for $f.$
The converse is, in general, false \cite{Frosini2} (see
also \cite{Frosini}).
\end{remark}
\begin{remark}\label{Nota su W_(G)cpa}
Let observe that $W_{G}(f)$ is arcwise connected. The proof is the same of proposition $3.14$ in \cite{Frosini2} (see
also \cite{Frosini}).
\end{remark}
In
order to state the  result  which characterizes the set $W_{G}(f), $
we need to introduce some definitions and results.
Herv\'e  proved the following
useful theorem (see
\cite{HerveSuBidisco} Theorem $1$):
\begin{theorem}\label{HerveTeo}
Let  $f:\Delta^{2}\to\Delta^{2}$ be a holomorphic map, without
 interior fixed points in $\Delta^{2}$,  whose components are
 $f_{1}, f_{2}.$ Then either
\begin{enumerate}
      \item  there exists a Wolff point,   $e^{\imath\theta_{1}}, $ of
             $f_{1}(\cdot\;, y)$,   which does  not depend on $y$
                                  or
      \item  there exists a holomorphic function
             $F_{1}:\Delta\rightarrow\Delta$,  such that

             $f_{1}(F_{1}(y), y)=F_{1}(y), \;\forall\;y\in\Delta.$
             In this case  $f_{1}(x, y)= x \Rightarrow$

             $x=F_{1}(y).$
\end{enumerate}
\end{theorem}
Let us remark that,  if $f\neq id_{\Delta^{2}}, $ then cases $i)$ and
$ii)$ cannot hold at the same time. Motivated by the
last mentioned  result of Herv\'e we give the following
de\-fi\-ni\-tion (\cite{Frosini}, \cite{Frosini2}):

\begin{definition}\label{MapType}
The holomorphic map $f:\Delta^{2}\to\Delta^{2},$ whose components are
 $f_{1}, f_{2},$  is called of:
\begin{enumerate}
\item \emph{first type} if:
      \begin{itemize}
       \item[-] there exists a holomorphic function
         $F_{1}:\Delta\rightarrow\Delta$,  such that

         $f_{1}(F_{1}(y), y)=F_{1}(y), \;\forall\;y\in\Delta$ and
       \item[-] there exists a holomorphic function
         $F_{2}:\Delta\rightarrow\Delta$,  such that

          $f_{2}(x, F_{2}(x))=F_{2}(x), \;\forall\;x\in\Delta.$
      \end{itemize}
\item \emph{second type} if (up to switching $f_{1}$ with $f_{2}$):
      \begin{itemize}
       \item[-] there exists a  Wolff point,  $e^{\imath\alpha_{1}}, $ of
                $f_{1}(\cdot\;, y)$,   (necessarily independent of $y$) and
       \item[-] there exists a holomorphic function
                $F_{2}:\Delta\rightarrow\Delta$,  such that

                $f_{2}(x, F_{2}(x))=F_{2}(x), \;\forall\;x\in\Delta.$
      \end{itemize}
\item \emph{third type} if:
      \begin{itemize}
       \item[-] there exists a Wolff point,  $e^{\imath\gamma_{1}}, $ of
                $f_{1}(\cdot\;, y)$,   (independent of $y$) and
       \item[-] there exists a Wolff point,   $e^{\imath\gamma_{2}}, $ of
                $f_{2}(x, \cdot)$,  (independent of $x$).
      \end{itemize}
\end{enumerate}
\end{definition}

In case $f$ is of \emph{first type} and without interior fixed points
in $\Delta^{2}$,  then it turns out that  $F_{1}\circ
F_{2}$ and $F_{2}\circ F_{1}$ have a Wolff point (see Lemma $3.10$ in \cite{Frosini2} and also
 \cite{Frosini}). Let $e^{\imath\theta_{1}} $
(respectively $e^{\imath\theta_{2}}$)  be the Wolff point of  $F_{1}\circ F_{2}$ (respectively $F_{2}\circ
F_{1}$). We also let $\lambda_{1}$ and $\lambda_{2}$  be,  respectively,   the \emph{boundary dilation coefficients} of $F_{1}$ at $
e^{\imath\theta_{2}}$ and of $F_{2}$ at $ e^{\imath\theta_{1}}$ (see Lemma$3.10$ in \cite{Frosini2} and also
 \cite{Frosini}). In case $f$ is of
\emph{second type} we denote  by
$e^{\imath\alpha_{1}}$ the Wolff point of $f_{1}(\cdot, y),$  by
$e^{\imath\alpha_{2}}$ the $K-$limit (or non-tangential limit) (see
definition in \cite{LibroAbate}) of $F_{2}$ at
$e^{\imath\alpha_{1}}$ (if it exists) and $k_{2}:=\lim\limits_{x\to e^{\imath\alpha_{1}}}|F_{2}'(x)|.$
In case $f$ is of \emph{third type}, we set $e^{\imath\gamma_{1}}$ and
$e^{\imath\gamma_{2}}$ to be, respectively, the Wolff points of
$f_{1}(\cdot, y)$ and $f_{2}(x, \cdot).$
We let $\pi_{j}:\Delta^{2}\to\Delta$ ($j=1, 2$) be the projection on the $j-$th component.
Finally without loss of generalization we suppose that $e^{\imath{\theta_{1}}}=e^{\imath{\theta_{2}}}
=e^{\imath{\alpha_{1}}}=e^{\imath{\alpha_{2}}}=e^{\imath\gamma_{1}}=e^{\imath\gamma_{2}}=1.$  With the
above established notations  we proved (see \cite{Frosini2}) the following result:

\begin{theorem}\label{W(f)}

Let $f=(f_{1},  f_{2}) $ be a holomorphic map,  without fixed points in the complex bidisc. If $f_{1}\neq\pi_{1}$
and $f_{2}\neq\pi_{2}, $ then only  the following five cases are possible:
\begin{itemize}
\item[i)] $W(f)=\emptyset$ if and only if $f$ is of
\emph{first type} and $\lambda_{i}>1$ for either $i=1$ or $i=2;$

\item[ii)] $W(f)={(1, 1)}$ iff $f$ is of \emph{first type} and
  $\lambda_{i}\leq 1$ for each  $i=1, 2;$

\item[iii)] $W(f)=\{\{1\} \times \Delta\}\cup \{(1,
1)\}$ iff $f$ is of
  \emph{second type} and
  $k_{2}\leq 1;$

\item[iv)]  $W(f)=\{\{1\} \times \Delta\}$ iff
$f$ is of \emph{second type} and
  $k_{2}> 1;$

\item[v)] $W(f)= \{\{1\}\times \Delta\}\cup \{(1,
1)\}\cup \{\Delta \times
  \{1\}\}$ iff $f$ is of \emph{third type}.

On the other hand,  if $f_{1}(x, y)=x, $ $\forall\; y\in \Delta,  $ i.e if $f_{1}=\pi_{1}$ (or respectively
$f_{2}(x, y)=y, $ $\forall\; x\in \Delta , $ i.e $f_{2}=\pi_{2}$ ) then:

\item[vi)] $W(f)=(1\times\Delta)\cup
(1,1) \cup (\Delta\times
1)\cup
(1,1)\cup
(1 \times\Delta)$  where $1$
is the Wolff point of $f_{2}(x,\cdot)$

(or respectively  $W(f)=(\Delta\times 1)\cup
(1,1) \cup
(1\times\Delta)\cup
(1,1)\cup (\Delta\times
1)$ where $1$ is the Wolff
point of $f_{1}(\cdot, y)$).
\end{itemize}
\end{theorem}
With the same techniques used in the proof of Theorem \ref{W(f)} (\cite{Frosini}, \cite{Frosini2}) we also get the following characterization
of the generalized Wolff points of $f$:
\begin{theorem}\label{W(f)generalized}

Let $f=(f_{1},  f_{2}) $ be a holomorphic map,  without fixed points in the complex bidisc. If $f_{1}\neq\pi_{1}$
and $f_{2}\neq\pi_{2}, $ then  the following three cases are possible:
\begin{itemize}
\item[i)] $W_{G}(f)={(1, 1)}$ if and only if $f$ is of
\emph{first type}.

\item[ii)]  $W_{G}(f)=\{\{1\} \times \Delta\}\cup \{(1,
1)\}$ iff $f$ is of
  \emph{second type};
\item[iii)] $W_{G}(f)= \{\{1\}\times \Delta\}\cup \{{(1,
1})\}\cup \{\Delta \times
  \{1\}\}$ iff $f$ is of \emph{third type}.

\end{itemize}

\end{theorem}

It is interesting to notice that, in this case, using the result of Herv\'{e} \cite{HerveSuBidisco},
about the target set of $f,$ we get that $W_{G}(f)\equiv T(f).$

\end{document}